\documentclass[1 [leqno,11pt]{amsart}
\usepackage{amssymb, amsmath}
\def\inte#1{
\displaystyle\mathop{#1\kern0pt}^\circ }

\def\ddj{\dot{\Delta}_j}
\def\dj{\Delta_j}

\def\s1{\sigma_1}

\let\d=\partial
\let\pa=\partial

\def\cB{{\mathcal B}}
\def\cC{{\mathcal C}}
\def\cD{{\mathcal D}}
\def\cE{{\mathcal E}}
\def\cF{{\mathcal F}}

\def\cH{{\mathcal H}}

\def\cS{{\mathcal S}}

\def\cV{{\mathcal V}}

\renewcommand{\div}{{\rm div}\,}

\newcommand{\loc}{{\rm loc}\,}

\newcommand{\Supp}{{\rm Supp}\,}

\newcommand{\Rmnum}[1]{\uppercase\expandafter{\romannumeral #1} }
 \numberwithin{equation}{section}

\def\dH{\dot{H}}
\def\dB{\dot{B}}

\def\virgp{\raise 2pt\hbox{,}}
\def\cdotpv{\raise 2pt\hbox{;}}

\def\eqdefa{\buildrel\hbox{\footnotesize def}\over =}

\def\C{\mathop{\mathbb C\kern 0pt}\nolimits}
\def\DD{\mathop{\mathbb D\kern 0pt}\nolimits}
\def\EE{\mathop{{\mathbb E \kern 0pt}}\nolimits}
\def\K{\mathop{\mathbb K\kern 0pt}\nolimits}
\def\N{\mathop{\mathbb N\kern 0pt}\nolimits}
\def\Q{\mathop{\mathbb Q\kern 0pt}\nolimits}
\def\R{\mathop{\mathbb R\kern 0pt}\nolimits}
\def\SS{\mathop{\mathbb S\kern 0pt}\nolimits}
\def\ZZ{\mathop{\mathbb Z\kern 0pt}\nolimits}
\def\TT{\mathop{\mathbb T\kern 0pt}\nolimits}
\def\PP{\mathop{\mathbb P\kern 0pt}\nolimits}

\newcommand{\Z}{{\ZZ}}

\def\dive{\mathop{\rm div}\nolimits}
\def\curl{\mathop{\rm curl}\nolimits}

\def\Supp{\mathop{\rm Supp}\nolimits\ }

\def\na{\nabla}
\def\pt{\partial_t}
\def\p3{\partial_3}
\def\ph{\partial_{\rm h}}
\def\Lh{\Delta_{\rm h}^{-1}}
\def\uh{u^{\rm h}}
\def\u3{u^3}
\def\bh{b^{\rm h}}
\def\b3{b^3}

\def\h{{\rm h}}
\def\v{{\rm v}}
\def\nh{\nabla_{\rm h}}

\def\xih{\xi_{\rm h}}

\def\om3{\omega^3}
\def\V{\rm V}
\def\omr{\omega_{\frac r2}}
\def\odr{(\Gamma_+)_{\frac r2}}
\def\omdr{(\Gamma_-)_{\frac r2}}

\def\htr{\cH^{\theta,r}}

\newcommand{\beq}{\begin{equation}}
\newcommand{\eeq}{\end{equation}}
\newcommand{\ben}{\begin{eqnarray}}
\newcommand{\een}{\end{eqnarray}}
\newcommand{\beno}{\begin{eqnarray*}}
\newcommand{\eeno}{\end{eqnarray*}}
\newcommand{\andf}{\quad\hbox{and}\quad}
\newcommand{\with}{\quad\hbox{with}\quad}
\newtheorem{defi}{Definition}[section]
\newtheorem{thm}{Theorem}[section]
\newtheorem{lem}{Lemma}[section]

\newtheorem{prop}{Proposition}[section]

\begin{document}
\title[Critical one component velocity regularity criteria to 3-D MHD system]
{On the critical one-component velocity regularity criteria to 3-D
incompressible MHD system}

\author[Y. Liu]{Yanlin Liu}
\address [Y. Liu]{department of mathematical sciences, university of science and technology of china, hefei 230026, china}
\email{liuyanlin3.14@126.com}


\date{\today}

\maketitle

\begin{abstract}
Let $(u,b)$ be a smooth enough solution of  3-D incompressible MHD
system.  We prove that if $(u,b)$ blows up at a finite time $T^*$,
then for any $p\in]4,\infty[$, there holds
$\int_0^{T^*}\bigl(\|u^3(t')\|^p_{\dH^{\frac 12+\frac
2p}}+\|b(t')\|^p_{\dH^{\frac 12+\frac 2p}}\bigr)dt'=\infty$. We
remark that all these quantities are in the critical regularity of
the MHD system.
\end{abstract}

\section{Introduction}\label{sec1}

In this work, we investigate necessary conditions for the breakdown
of regular solutions to the following 3-D incompressible
Magnetohydrodynamics (MHD in short) system
\begin{equation}\label{MHD}
\left\{
\begin{array}{l}
\displaystyle \d_t u+u\cdot\nabla u-b\cdot\nabla b+\nabla p=\Delta u,   \qquad(t,x)\in\R^+\times\R^3\\
\displaystyle \d_t b+u\cdot\nabla b-b\cdot\nabla u=\Delta b,\\
\displaystyle \div u=\div b=0,\\
\displaystyle u|_{t=0} =u_0,b|_{t=0} =b_0,
\end{array}
\right.
\end{equation}
where u, p denote the velocity and scalar pressure of the fluid
respectively, and b denotes the magnetic field.

\smallbreak When the initial magnetic field $b_0$ is identically
zero, the system \eqref{MHD} reduces to the classical Navier-Stokes
equations, the global regularity of which is still one of the
biggest open questions in the field of mathematical fluid mechanics.
Of course, the analogous problem for the MHD system remains just as
difficult due to the coupling with the magnetic field.

\medbreak This system has two major basic  features. First of all,
 the
total kinetic energy is conserved for smooth enough solutions of
\eqref{MHD}  \beq \label {kineticenergy} \begin{split} \frac 12
\bigl(\|u(t)\|_{L^2}^2+\|b(t)\|^2_{L^2}\bigr)
+\int_{0}^t\bigl(\|\nabla u(t')\|_{L^2}^2&+\|\na
b(t')\|_{L^2}^2\bigr) dt'=\frac 12
\bigl(\|u_0\|_{L^2}^2+\|b_0\|^2_{L^2}\bigr).
\end{split} \eeq The second basic feature is the scaling invariance.
 Indeed, if~$(
u,b,p)$ is a solution of \eqref{MHD} on~$[0,T]\times \R^3$,
then~$(u,b,p)_\lambda$ defined by \beq \label{scaling}
 (u,b,p)_\lambda(t,x) \eqdefa \bigl(
\lambda u (\lambda^2t, \lambda x), \lambda b (\lambda^2t, \lambda
x), \lambda^2 p(\lambda^2t,\lambda x)\bigr) \eeq is also a solution
of \eqref{MHD} on $[0,\lambda^{-2} T]\times \R^3$. This leads to the
notion of critical regularity corresponding to the System
\eqref{MHD}.

\smallbreak Before Proceeding,  let us  set
\begin{equation}\label{S1eq1}
\Omega\eqdefa\nabla\times u,\ \ j\eqdefa\nabla\times b,\ \
\omega\eqdefa\Omega\cdot e^3,\ \ d\eqdefa j\cdot e^3\  \
\mbox{with}\ \ e^3=(0,0,1).
\end{equation}

Motivated by the critical one component criteria in \cite{CZ14} by
Chemin and Zhang for the 3-D classical Navier-Stokes system,
Yamazaki \cite{Yamazaki} proved the following regularity criteria
for the System \eqref{MHD}:

\begin{thm}\label{thm1.1}
{\sl Let $\Omega_0,j_0\in L^\frac 32(R^3)$. Then the MHD system
\eqref{MHD} has a unique  solution $(u,b)$ on $[0,T^*[$  such that
$u,b\in C([0,T^*[;\dH^\frac 12(R^3))\bigcap
L^2_{\loc}(]0,T^*[;\dH^\frac 32(R^3))$ and
\begin{align}\label{1.2}
\sup_{t\in[0,T]}\bigl(\|\Omega(t)\|^{\frac 32}_{L^\frac
32}+\|j(t)\|^{\frac 32}_{L^\frac
32}\bigr)+\int_0^T\int_{R^3}\bigl(|\nabla(\Omega+j)|^2|\Omega+j|^{-\frac
12}+
|\nabla(\Omega-j)|^2|\Omega-j|^{-\frac 12}\bigr)dx dt'\notag\\
\leqslant C\bigl(1+\|\Omega_0\|^{\frac 32}_{L^\frac
32}+\|j_0\|^{\frac 32}_{L^\frac 32}\bigr)
\exp\Bigl(\int_0^T\bigl(\|u(t')\|^2_{\dH^{\frac
32}}+\|b(t')\|^2_{\dH^{\frac 32}}\bigr)dt'\Bigr)<\infty,
\end{align} for any $ T<T^*.$
Moreover, for $p\in]4,6[,p_1>9,p_2>\frac 92$, we denote
$$\|b\|_{SC_{p,p_1,p_2}}\eqdefa\|b\|^p_{\dH^{\frac 12+\frac 2p}}
+\|b\|^{r_1}_{L^{p_1}}+\|\nabla b\|^{r_2}_{L^{p_2}},\with
\frac{3}{p_1}+\frac{2}{r_1}=1,\frac{3}{p_2}+\frac{2}{r_2}=2.$$ If
$T^*<\infty$, then
\begin{equation}\label{1.3}
\int_0^{T^*}\bigl(\|u^3(t')\|^p_{\dH^{\frac 12+\frac
2p}}+\|b(t')\|_{SC_{p,p_1,p_2}}\bigr)dt'=\infty.
\end{equation}}
\end{thm}
It is easy to check that when $T^*=\infty,$ the quantity \eqref{1.3}
is scaling invariant under the scaling transformation
\eqref{scaling}.

\smallbreak The main result in \cite{CZ14} states that if $u$ is a
Fujta-Kato type solution to the classical Navier-Stokes system on
$[0,T^*[$ and if $T^*<\infty,$ then \eqref{1.3} holds for $p\in
]4,6[$ with $b=0.$ Very recently, this result was extended by
Chemin, Zhang and Zhang in \cite{CZZ} for $p\in ]4,\infty[.$
Corresponding to \cite{CZZ}, the purpose of this work is to extend
$p$ in Theorem \ref{thm1.1} to be in $]4,\infty[$ and to get rid of
the terms $\|b\|^{r_1}_{L^{p_1}}+\|\nabla b\|^{r_2}_{L^{p_2}}$ in
\eqref{1.3} by using the symmetric structure of the MHD system
\eqref{MHD}. One may check \cite{Yamazaki} and the references
therein for the other types of regularity criteria for the MHD
system (see \cite{CW10} for instance).

\smallbreak In all that follows, we consider initial data
$(u_0,b_0)$ with $\Omega_0,j_0\in L^\frac 32(R^3)$ so that Theorem
\ref{thm1.1} always holds. We shall concentrate on the proof of the
extended regularity criterion. In order to do so, let us recall the
following family of spaces from \cite{CZZ}.

\begin{defi}\label{def1.1}
For r in $[3/2, 2]$, we denote by $\cV^r$ the space of divergence
free vector fields with the vorticity of which belongs to $L^{\frac
32}\bigcap L^r$.
\end{defi}

Let us remark that, if we denote
\begin{equation}\label{1.4}
\alpha(r)\eqdefa\frac 1r -\frac 12,
\end{equation}
the dual Sobolev embedding $L^r\hookrightarrow\dH^{-3\alpha (r)}$
together with Biot-Savart's law implies that  $\cV^r$  is a subspace
of $\dH^{\frac 12}\bigcap\dH^{1-3\alpha(r)}$.

Our main result states as follows:

\begin{thm}\label{thm1.2}
Let us consider initial data $u_0,b_0\in\cV^2$. If the lifespan
$T^*$ of the unique maximal solution  (u,b) given by Theorem
\ref{thm1.1} is finite, then for any $p\in]4,\infty[$, we have
\begin{equation}\label{blouupc}
\int_0^{T^*}\left(\|u^3(t')\|^p_{\dH^{\frac 12+\frac
2p}}+\|b(t')\|^p_{\dH^{\frac 12+\frac 2p}}\right)dt'=\infty.
\end{equation}
\end{thm}

Let us complete this introduction by the notations we shall use in
the whole text.

Let $A, B$ be two operators, we denote $[A;B]=AB-BA,$ the commutator
between $A$ and $B$. For $a\lesssim b$, we mean that there is a
uniform constant $C,$ which may be different on different lines,
such that  $a\leq Cb,$ and $a\sim b$ means that both $a\lesssim b$
and $b\lesssim a$ hold. C stands for some universal positive constant
which may change from line to line
 and $C_0$ denotes a positive constant depending on the initial data
only. For a Banach space B, we shall use the shorthand $L^p_t(B)$
for $L^p(]0,t[;B)$.

\setcounter{equation}{0}
\section{Scheme of the proof and the organization of the paper.}

In fact, we shall prove the following more general version of
Theorem \ref{thm1.2}:

\begin{thm}\label{thm1.3}
{\sl Let $r\in[3/2,2[$ and $u_0,b_0\in\cV^r$. If the lifespan $T^*$
of the unique maximal solution $(u,b)$  given by Theorem
\ref{thm1.1} is finite, then for any $p\in
\bigl]4,\frac{2r}{2-r}\bigr[$, we have
\begin{equation}\label{1.5}
\int_0^{T^*}\left(\|u^3(t')\|^p_{\dH^{\frac 12+\frac
2p}}+\|b(t')\|^p_{\dH^{\frac 12+\frac 2p}}\right)dt'=\infty.
\end{equation}
}
\end{thm}

The main idea of the proof here basically follow from \cite{CZ14,
CZZ, Yamazaki}. We first recall some important definitions and
notations. Let
$$\nabla_\h^\bot\eqdefa(-\partial_2,\partial_1),\andf \Delta_\h\eqdefa\partial_1^2+\partial_2^2,$$
for any $f=(f^\h,f^3)$ with $\dive f=0$, we write
\begin{equation}\label{1.6}
f^\h=f^\h_{\curl}+f^\h_{\div}\ \with \
f^\h_{\curl}\eqdefa\nh^\bot\Lh(\nabla\times f)\cdot e^3,\ \
f^\h_{\div}\eqdefa-\nh\Lh\partial_3 f^3.
\end{equation}
This is sort of Hodge decomposition for the horizontal variables,
and we emphasize that this is a key identity to be used frequently
in what follows. Moreover, because of the operator $\nh\Lh$, it is
naturally to measure horizontal derivatives and vertical derivatives
differently. This leads to the following definition of the
anisotropic Sobolev spaces.

\begin{defi}[Definition 2.1 of \cite{CZ14,CZZ}]\label{def1.2}
For   $(s,s')$ in $\R^2$, $\dH^{s,s^\prime}$ denotes the space of
tempered distribution a such that
$$\|a\|_{\dH^{s,s'}}^2\eqdefa\int_{\R^3}|\xi_\h|^{2s}|\xi_3|^{2s'}|\widehat{a}(\xi)|^2d\xi<\infty\ with\ \xi_\h=(\xi_1,\xi_2).$$
For $\alpha(r)$ given by \eqref{1.4} and $\theta\in ]0,3\alpha(r)[$,
we denote $\htr\eqdefa\dH^{-3\alpha(r)+\theta,-\theta}.$
\end{defi}
Then it follows from (2.7) of \cite{CZZ} that
\begin{equation}\label{1.7}
\|\p3 u^3\|_{\htr}\lesssim\|u\|_{\dH^{1-3\alpha(r)}}.
\end{equation}
To use the space efficiently in the proof, we need to rely them on
anisotropic Littlewood-Paley theory and also anisotropic Besov
spaces. These will be done in the following section.

The first step to prove Theorem \ref{thm1.3} is the following
proposition:

\begin{prop}\label{prop1.1}
{\sl Under the hypothesis of Theorem \ref{thm1.3}, for any $p\in
\bigl]4,\frac{2r}{2-r}\bigr[,\, \theta\in ]0,\alpha(r)[$, a constant
$C$ exists such that for any $t<T^*$, we have \beq\label{1.8}
\begin{split}
\frac 1r\bigl(\|& (\Gamma_+)_{\frac{r}2}(t)\|_{L^2}^2+\|\omdr(t)\|_{L^2}^2\bigr)\\
&+\frac{2(r-1)}{r^2}\int_0^t\bigl(\|\nabla\odr(t')\|_{L^2}^2+\|\nabla \omdr(t')\|_{L^2}^2\bigr)dt'\\
 \lesssim &\Bigl(\frac 2r\bigl(\|\omega_0\|_{L^r}^r+\|d_0\|_{L^r}^r\bigr)\\
&\quad+\Bigl(\int_0^t\|\nabla\p3{\V_+}(t')\|_{\htr}^2
+\|\nabla\p3{\V_-}(t')\|_{\htr}^2 dt'\Bigr)^{\frac
r2}\Bigr)\cdot\cE(t).
\end{split} \eeq
Here and in all that follows, we always denote \beq \label{1.10}
\begin{split} & \Gamma_+\eqdefa \omega+d,\quad \Gamma_-\eqdefa \omega-d,\quad
\V_+\eqdefa u^3+b^3,\quad \V_-\eqdefa u^3-b^3 \andf
\\
& a_{\alpha}\eqdefa\frac{a}{|a|}|a|^{\alpha},\quad
\cE(t)\eqdefa\exp\Bigl(C\int_0^t\bigl(\|u^3(t')\|_{\dH^{\frac
12+\frac 2p}}^p+\|b(t')\|_{\dH^{\frac 12+\frac
2p}}^p\bigr)dt'\Bigr). \end{split} \eeq
  for
scalar function a and $\alpha\in]0,1[$ } \end{prop}

To prove this proposition,  we need to use the structures of the
equations for $\omega$ and d, namely  \eqref{3.1}.  The quadratic
terms $\uh_{\curl}\cdot\nh\omega$ and $\uh_{\curl}\cdot\nh d$ look
dangerous. As in \cite{CZ14, CZZ, Yamazaki}, a way to get rid of it
is to use an energy type estimate and the divergence-free condition.
Here we shall perform an $L^r$ energy estimate for $\omega$ and d based on the following lemma.

\begin{lem}[Lemma 3.1 of \cite{CZ14}]\label{lem1.1}
{\sl  Let r be in $]1,2[$ and $a_0$ a function in $L^r$. Let us
consider a function f in $L_{loc}^1(\R^+;L^r)$ and $v$ a divergence
free vector field in $L_{loc}^2(\R^+;L^{\infty})$. If a solves
\begin{equation*}
\left\{
\begin{array}{c}
\pt a-\Delta a+v\cdot\nabla a=f\\
a|_{t=0}=a_0,
\end{array}
\right.
\end{equation*}
then $|a|^{\frac r2}$ belongs to $L_{loc}^{\infty}(\R^+;L^2)\bigcap
L_{loc}^2(\R^+;\dH^1)$ and
\begin{align*}
\frac 1r\int_{\R^3}|a(t,x)|^r& dx+(r-1)\int_0^t\int_{\R^3}|\nabla a(t',x)|^2|a(t',x)|^{r-2} dxdt'\\
& =\frac 1r\int_{\R^3}|a_0(x)|^r
dx+\int_0^t\int_{\R^3}f(t',x)a(t',x)|a(t',x)|^{r-2} dxdt'.
\end{align*}}
\end{lem}
The proof of Proposition \ref{prop1.1} is the purpose of the
fourth section.

We remark that for the MHD system \eqref{MHD}, additional difficulty
arises  in the estimate of $\|\nabla\p3
{\V_+}\|_{L_t^2(\htr)}+\|\nabla \p3 {\V_-}\|_{L_t^2(\htr)}$ due to
the appearance of terms like
$2(\pa_1\bh\cdot\pa_2\uh-\pa_2\bh\cdot\pa_1\uh)$ in right-hand side
of  \eqref{3.1}. This is the purpose of the next proposition.

\begin{prop}\label{prop1.2}
{\sl Under the hypothesis of Theorem \ref{thm1.3}, for any
$p\in\bigl]4,\frac{2r}{2-r}\bigr[,\, \theta\in\bigl]3\alpha(r)-\frac
2p,\alpha(r)\bigr[$, a constant C exists such that for any $t<T^*$, we
have \beq\label{1.11}
\begin{split}
\|\p3 &{\V_+}(t) \|_{\htr}^2+\|\p3 {\V_-}(t)\|_{\htr}^2+\int_0^t\bigl(\|\nabla\p3 {\V_+}(t')\|_{\htr}^2+\|\nabla \p3 {\V_-}(t')\|_{\htr}^2\bigr)dt'\\
\lesssim
&\left(\|\Omega_0\|_{L^r}^2+\|j_0\|_{L^r}^2+\int_0^t\Bigl(\bigl(\|\u3\|_{\dH^{\frac
12+\frac 2p}}
+\|\b3\|_{\dH^{\frac 12+\frac 2p}}\bigr)\right.\\
&\times\bigl(\|(\Gamma_+)_{\frac r2} \|_{L^2}^{2(2\alpha(r)+\frac 1p)}\|\nabla\odr\|_{L^2}^{2(1-\frac 1p)}+\|\omdr\|_{L^2}^{2(2\alpha(r)+\frac 1p)}\|\nabla \omdr\|_{L^2}^{2(1-\frac 1p)}\bigr)\\
&+\bigl(\|u^3\|_{\dH^{\frac 12+\frac 2p}}^2+\|b^3\|_{\dH^{\frac 12+\frac 2p}}^2\bigr)\bigl(\|\odr\|_{L^2}^2+\|\omdr\|_{L^2}^2)^{2\bigl(\alpha(r)+\frac 1p\bigr)}\\
&\left.\times \bigl(\|\nabla\odr\|_{L^2}^2+\|\nabla
\omdr\|_{L^2}^2)^{1-\frac 2p}\Bigr)dt'\right)\cE(t).
\end{split} \eeq
}
\end{prop}
The proof of Proposition \ref{prop1.2} is the purpose of the fifth section.

Finally we close the estimates by the following proposition:
\begin{prop}\label{prop1.3}
{\sl Under the hypothesis of Theorem \ref{thm1.3}, for any $p\in
\bigl]4,\frac{2r}{2-r}\bigr[,\, \theta\in \bigl]3\alpha(r)-\frac
2p,\alpha(r)\bigr[$, a constant C exists such that for any $t<T^*$,
we have
\begin{align}\label{1.12}
\begin{split}
&\|\odr(t) \|_{L^2}^{2(1+2p\alpha(r))}+\|\omdr(t)\|_{L^2}^{2(1+2p\alpha(r))}+\|\nabla\odr\|_{L_t^2L^2}^{2(1+2p\alpha(r))}\\
&\qquad\quad+\|\nabla \omdr\|_{L_t^2L^2}^{2(1+2p\alpha(r))} \lesssim
\bigl(\|\Omega_0\|_{L^r}^{r(1+2p\alpha(r))}+\|j_0\|_{L^r}^{r(1+2p\alpha(r))}\bigr)\exp\bigl(C\cE(t)\bigr),
\end{split}
\end{align}
and
\begin{align}\label{1.13}
\begin{split}
\|\p3& {\V_+}(t)\|_{\htr}^2+\|\p3 {\V_-}(t)\|_{\htr}^2\\
&+\int_0^t\bigl(\|\nabla \p3 {\V_+}(t')\|_{\htr}^2+\|\nabla
\p3{\V_-}(t')\|_{\htr}^2\bigr)dt\\
\lesssim&\bigl(\|\Omega_0\|_{L^r}^2+\|j_0\|_{L^r}^2\bigr)\exp\bigl(C\cE(t)\bigr).
\end{split}
\end{align}}
\end{prop}
The proof of Proposition \ref{prop1.3} is the purpose of the sixth section.

Now we have controls on the quantities \beq\label{6.2}
\begin{split}
& \sup_{t\in [0,T^\ast[}\|(\Gamma_+)(t)\|_{L^r},\ \int_0^{T^*}\|\nabla\odr\|_{L^2}^2dt',\ \int_0^{T^*}\|\nabla\p3{\V_+}\|_{\htr}^2dt'\andf\\
& \sup_{t\in [0,T^\ast[}\|(\Gamma_-)(t)\|_{L^r},\
\int_0^{T^*}\|\nabla\omdr\|_{L^2}^2dt',\
\int_0^{T^*}\|\nabla\p3{\V_-}\|_{\htr}^2dt'.
\end{split}
\eeq We want to prove that all the above quantities prevent the
solution of \eqref{MHD} from blowing up. The details will be
presented in the last section.

\section{Preliminaries}\label{sec2}

In this section, we  first recall some basic facts on anisotropic
Littlewood-Paley theory from \cite{BCD, cz07, Pa05}, and then we
collect some interesting estimates from \cite{CZ14, CZZ} that will
be used later on.

\subsection{Basic facts on Littlewood-Paley theory}\label{sec2.1}
 Let
$\cC\eqdefa\{\xi\in\R^d:\frac 34\leqslant|\xi|\leqslant\frac 83\}$
and $\cB\eqdefa \{\xi\in\R^d:|\xi|\leqslant\frac 43\}$. There exist
two radial functions $\chi\in\cD(\cB)$ and $\varphi\in\cD(\cC)$ such
that
$$\chi(\xi)+\sum_{j\geqslant0}\varphi(2^{-j}\xi)=1,\ \forall\ \xi\in\R^d\quad and\quad \sum_{j\in\Z}\varphi(2^{-j}\xi)=1,\ \forall\ \xi\in\R^d\backslash\{0\}.$$
For every $a\in\cS'(\R^3),$ we recall the dyadic operator for both
isentropic and anisotropic version
\begin{align}\label{2.1}
&\ddj a=\cF^{-1}(\varphi(2^{-j}|\xi|)\widehat{a}),\qquad \  \ \dot{S}_j a=\cF^{-1}(\chi(2^{-j}|\xi|)\widehat{a}),\notag\\
&\dot{\Delta}^\h_k a=\cF^{-1}(\varphi(2^{-k}|\xih|)\widehat{a}),\qquad \dot{S}_k^\h a=\cF^{-1}(\chi(2^{-k}|\xih|)\widehat{a}),\\
&\dot{\Delta}^\v_\ell
a=\cF^{-1}(\varphi(2^{-\ell}|\xi_3|)\widehat{a}),\qquad \
\dot{S}_\ell^\v a=\cF^{-1}(\chi(2^{-\ell}|\xi_3|)\widehat{a}),\notag
\end{align}
where $\xih=(\xi_1,\xi_2)$, $\cF a$ and $\widehat{a}$ denote the
Fourier transform of a.

 Moreover, it is easy to verify that
for any u in $\cS'_h$, which means that $u$ belongs to  $\cS'$ and
satisfies $\lim\limits_{j\rightarrow-\infty}\|\dot{S}_ju\|_{L^\infty}=0,$
there holds $u=\sum\limits_{j\in\mathbb{Z}}\ddj u.$

 Let us recall the homogeneous isentropic Besov space from \cite{BCD}.

\begin{defi}\label{def2.1}
Let $1\leqslant p,r\leqslant+\infty$ and $s\in\R$. For any u in
$\cS'_h(\R^3)$, we set
$$\|u\|_{\dB^s_{p,r}}\eqdefa\|(2^{js}\|\ddj u\|_{L^p})_j\|_{\ell^r(\mathbb{Z})}.$$
\begin{itemize}
  \item For $s<\frac 3p$ (or $s=\frac 3p$ if r=1), we define $\dB^s_{p,r}(\R^3)\eqdefa\{u\in\cS'_h(\R^3)\,:\,\|u\|_{\dB^s_{p,r}}<\infty\}$.
  \item If there exists some positive integer k such that $\frac 3p+k\leqslant s<\frac 3p+k+1$ (or $s=\frac 3p+k+1$ if r=1), then we define $\dB^s_{p,r}(\R^3)$ as the subset of distributions u in $\cS'_h(\R^3)$ such that $\pa^{\beta}u$ belongs to $\dB^{s-k-1}_{p,r}$ whenever $|\beta|=k+1$.
\end{itemize}
\end{defi}
We remark that in particular, $\dB^s_{2,2}$
coincides with the classical homogeneous Sobolev space $\dH^s$.
Similarly, we can also define the homogeneous anisotropic Besov
space.
\begin{defi}\label{def2.2}
Let us define the homogeneous anisotropic Besov space
$(\dB^{s_1}_{p,q_1})_\h(\dB^{s_2}_{p,q_2})_\v$ as the subspace of
distributions u in $\cS'_h(\R^3)$ such that
$$\|u\|_{(\dB^{s_1}_{p,q_1})_\h(\dB^{s_2}_{p,q_2})_\v}\eqdefa
\Bigl(\sum_{k\in\mathbb{Z}}2^{q_1ks_1}\bigl(\sum_{\ell\in\mathbb{Z}}2^{q_2\ell
s_2}\|\dot{\Delta}_k^\h\dot{\Delta}_\ell^\v
u\|_{L^p}^{q_2}\bigr)^{\frac{q_1}{q_2}}\Bigr)^{\frac{1}{q_1}}$$ is
finite.
\end{defi}
We remark that in particular,
$(\dB^{s_1}_{2,2})_\h(\dB^{s_2}_{2,2})_\v$ coincides with the
homogeneous anisotropic Sobolev sapce $\dH^{s_1,s_2}$, and thus the
space $(\dB^{-3\alpha(r)+\theta}_{2,2})_\h(\dB^{-\theta}_{2,2})_\v$
is the space $\htr$ given by Definition \ref{def1.2}. Let us also
remark that in the case when $q_1$ is different from $q_2$, the
order of summation is important.

By virtue of the above definitions,  one has

\begin{lem}[Lemma 4.3 of \cite{CZ14}]\label{lem2.1}
{\sl For any $s>0$ and any $\theta\in ]0,s[$, we have
\begin{equation}\label{2.2}
\|f\|_{(\dB^{s-\theta}_{p,q})_\h(\dB^{\theta}_{p,1})_\v}\lesssim\|f\|_{\dB^s_{p,q}}.
\end{equation}}
\end{lem}

We also recall the following Bernstein type lemmas:
\begin{lem}[Isentropic version, see \cite{BCD}]\label{lem2.2}
{\sl  Let $\cC$ be an annulus and $\cB$ a ball of $\R^3$. Then for
any nonnegative integer N, and $1\leqslant p\leqslant
q\leqslant\infty$, we have
$$\Supp \widehat{a}\subset\lambda\cB\Longrightarrow\|D^N a\|_{L^q}\eqdefa\sup_{|\alpha|=N}\|\pa^\alpha a\|_{L^q}
\lesssim\lambda^{N+3(\frac 1p-\frac 1q)}\|a\|_{L^p},$$
$$\Supp \widehat{a}\subset\lambda\cC\Longrightarrow \lambda^N\|a\|_{L^p}\lesssim\|D^N
a\|_{L^p}\lesssim\lambda^N\|a\|_{L^p}.$$}
\end{lem}

\begin{lem}[Anisotropic version, see \cite{cz07, Pa05}]\label{lem2.3}
{\sl  Let $\cC_\h$ (resp. $\cC_\v$) be an annulus of $\R^2_\h$
(resp. $\R_\v$), and $\cB_\h$ (resp. $\cB_\v$) a ball of $\R^2_\h$
(resp. $\R_\v$). Then for any nonnegative integer N, and $1\leqslant
p_2\leqslant p_1\leqslant\infty$ and $1\leqslant q_2\leqslant
q_1\leqslant\infty$, we have \beno
\begin{split}
\Supp \widehat{a}\subset\lambda\cB_\h&\Longrightarrow\|\ph^\alpha
a\|_{L_\h^{p_1}(L_\v^{q_1})}\lesssim\lambda^{|\alpha|+2(\frac
{1}{p_2}-\frac {1}{p_1})}\|a\|_{L_\h^{p_2}(L_\v^{q_1})},\\
\Supp \widehat{a}\subset\lambda\cB_\v&\Longrightarrow\|\p3^\beta
a\|_{L_\h^{p_1}(L_\v^{q_1})}\lesssim\lambda^{|\beta|+(\frac
{1}{q_2}-\frac {1}{q_1})}\|a\|_{L_\h^{p_1}(L_\v^{q_2})},\\
\Supp
\widehat{a}\subset\lambda\cC_\h&\Longrightarrow\|a\|_{L_\h^{p_1}(L_\v^{q_1})}
\lesssim\lambda^{-N}\sup_{|\alpha|=N}\|\ph^{\alpha}a\|_{L_\h^{p_1}(L_\v^{q_1})},\\
\Supp
\widehat{a}\subset\lambda\cC_\v&\Longrightarrow\|a\|_{L_\h^{p_1}(L_\v^{q_1})}
\lesssim\lambda^{-N}\|\p3^N a\|_{L_\h^{p_1}(L_\v^{q_1})}.\end{split}
\eeno}
\end{lem}
As a corollary of Lemma \ref{lem2.3}, for any $1\leqslant
p_2\leqslant p_1\leqslant\infty$, we have
$$\|a\|_{(\dB^{s_1-2(\frac{1}{p_2}-\frac{1}{p_1})}_{p_1,q_1})_\h(\dB^{s_2-(\frac{1}{p_2}-\frac{1}{p_1})}_{p_1,q_2})_\v}
\lesssim\|a\|_{(\dB^{s_1}_{p_2,q_1})_\h(\dB^{s_2}_{p_2,q_2})_\v}.$$


\subsection{Some technical inequalities}\label{sec2.2} For the
convenience of the readers, we recall some inequalities from
\cite{{CZ14},{CZZ}} that will be used in what follows.

\begin{lem}[Lemma 3.1 of \cite{CZZ}]\label{lem2.4}
{\sl For r in $]3/2,2[$, we have
\begin{equation}\label{2.3}
\|\nabla a\|_{L^r}\lesssim\|\nabla a_{\frac r2}\|_{L^2}\|a_{\frac
r2}\|_{L^2}^{\frac 2r-1}.
\end{equation}
Moreover, for s in $[-3\alpha(r),1-\alpha(r)]$, we have
\begin{equation}\label{2.4}
\|a\|_{\dH^s}\lesssim\|\nabla a_{\frac
r2}\|_{L^2}^{3\alpha(r)+s}\|a_{\frac r2}\|_{L^2}^{1-\alpha(r)-s}.
\end{equation}}
\end{lem}

\begin{lem}[Proposition 3.1 of \cite{CZZ}]\label{lem2.5}
{\sl Let $u$ be a divergence-free vector field. For
$\theta\in]0,3\alpha(r)[$ and $\beta\in]0,1/2[$, we have
\begin{equation}\label{2.5}
\|\uh\|_{(\dB^1_{2,1})_\h(\dB^{1-3\alpha(r)-\beta}_{2,1})_\v}\lesssim\|\omr\|_{L^2}^{2\alpha(r)+\beta}\|\nabla\omr\|_{L^2}^{1-\beta}
+\|\p3\u3\|_{\htr}^{\beta}\|\nabla\p3\u3\|_{\htr}^{1-\beta}.
\end{equation}}
\end{lem}

It is easy to observe that the proof of Lemma 5.2 in \cite{CZZ}
implies the following inequality:

\begin{lem}\label{lem2.6}
{\sl Let $f=(f^1,f^2,f^3),g=(g^1,g^2,g^3)$ and $f^\h=(f^1,f^2)$.
Then for any p in $\bigl]4,\frac{2r}{2-r}\bigr[$ and any $\theta$ in
$]3\alpha(r)-2/p,\alpha(r)[$, we have
\begin{align}\label{2.6}
\begin{split}
|(f^\h\cdot\nh\p3 g^3|\p3 g^3)& _{\htr}|\lesssim\|g^3\|_{\dH^{\frac 12+\frac 2p}}\Bigl(\|\nh f^\h\|^2_{\dH^{\frac 12-3\alpha(r)+\theta,\frac 12-\frac 1p-\theta}}\\
& +\|\p3 g^3\|^2_{\dH^{\frac 12-3\alpha(r)+\theta,\frac 12-\frac
1p-\theta}} +\|f^\h\|_{(\dB^1_{2,1})_\h(\dB^{1-3\alpha(r)-\frac
2p}_{2,1})_\v}\|\nabla\p3 g^3\|_{\htr}\Bigr).
\end{split}
\end{align}}
\end{lem}

\section{Proof of Proposition \ref{prop1.1}}\label{sec3}
The purpose of this section is to present the proof of Proposition
\ref{prop1.1}. Note that $\omega=\pa_1 u^2-\pa_2 u^1,\ d=\pa_1
b^2-\pa_2 b^1$, then it follows from \eqref{MHD} that
\begin{align*}
& \pt\omega+u\cdot\nabla\omega-b\cdot\nabla d-\Delta\omega=\pa_1
u\cdot\nabla u^2-\pa_2 u\cdot\nabla u^1-\pa_1 b\cdot\nabla b^2
+\pa_2 b\cdot\nabla b^1,\\
& \pt d +u\cdot\nabla d-b\cdot\nabla\omega-\Delta d=\pa_1
u\cdot\nabla b^2-\pa_2 u\cdot\nabla b^1-\pa_1 b\cdot\nabla u^2
+\pa_2 b\cdot\nabla u^1,
\end{align*}
from which, and $\dive u=\dive b=0,$ we deduce
\begin{equation}\label{3.1}
\left\{
\begin{array}{l}
\pt\omega+u\cdot\nabla\omega-b\cdot\nabla d-\Delta\omega=(\p3 u^3\omega+\pa_2 u^3\p3 u^1-\pa_1 u^3\p3 u^2)\\
\qquad-(\p3 b^3 d+\pa_2 b^3\p3 b^1-\pa_1 b^3\p3 b^2),\\
\pt d +u\cdot\nabla d-b\cdot\nabla\omega-\Delta d=\p3\u3 d-\p3\b3\omega-\pa_1\u3\p3 b^2+\pa_2\u3\p3 b^1\\
\qquad+\pa_1\b3\p3 u^2-\pa_2\b3\p3
u^1+2(\pa_1\bh\cdot\pa_2\uh-\pa_2\bh\cdot\pa_1\uh).
\end{array}
\right.
\end{equation}
Summing up these two equations  gives \beq\label{3.2}\begin{split}
&\pt{\Gamma_+}+u\cdot\nabla{\Gamma_+}-b\cdot\nabla{\Gamma_+}-\Delta{\Gamma_+}
=\p3 \V_-{\Gamma_+}\\
 &\qquad\qquad\quad-\pa_1\V_-\p3
(u^2+b^2)+\pa_2\V_-\p3
(u^1+b^1)+2(\pa_1\bh\cdot\pa_2\uh-\pa_2\bh\cdot\pa_1\uh).
\end{split}\eeq

Since $\dive u=\dive b=0,$ we get,  by applying Lemma \ref{lem1.1},
that \begin{equation}\label{3.3}
\frac
1r\|\odr(t)\|_{L^2}^2+\frac{4(r-1)}{r^2}\int_0^t\|\nabla\odr\|_{L^2}^2dt'=\frac
1r\bigl\|(\omega_0+d_0)\bigr\|_{L^r}^r+\sum_{i=1}^3 I_i,
\end{equation}
where
\begin{align}
&I_1=\int_0^t\int\p3{\V_-}{\Gamma_+}{(\Gamma_+)}_{r-1} dxdt',\notag\\
&I_2=\int_0^t\int\bigl(-\pa_1{\V_-}\p3(u^2+b^2)+\pa_2{\V_-}\p3(u^1+b^1)\bigr){(\Gamma_+)}_{r-1} dxdt',\\
&I_3=2\int_0^t\int(\pa_1\bh\cdot\pa_2\uh-\pa_2\bh\cdot\pa_1\uh){(\Gamma_+)}_{r-1}
dxdt'.\notag
\end{align}
We first get, by using integrating by parts, that
\begin{align*}
|I_1|& \leqslant r\int_0^t\int|{\V_-}||\p3{\Gamma_+}||{\Gamma_+}|^{r-1} dxdt'\\
& =r\int_0^t\int|{\V_-}||\p3{\Gamma_+}||\odr|^{\frac{2}{r'}}dxdt'\\
& \leqslant
r\int_0^t(\|\u3\|_{L^{\frac{3p}{p-2}}}+\|\b3\|_{L^{\frac{3p}{p-2}}})\|\p3{\Gamma_+}\|_{L^r}
\|\odr\|_{L^{\frac{6p(r-1)}{2pr-3p+2r}}}^{\frac{2}{r'}}dt',
\end{align*}
where $r'$ denotes the conjugate index of r so that $\frac
1r+\frac{1}{r'}=1$. As $p\in \bigl]4,\frac{2r}{2-r}\bigr[$, we have
that $r'\cdot\frac{p-2}{2p}\in ]0,1[$, then Sobolev embedding and
interpolation inequality imply that
$$\|\odr\|_{L^{\frac{6p(r-1)}{2pr-3p+2r}}}\lesssim\|\odr\|_{\dH^{r'\cdot\frac{p-2}{2p}}}\lesssim\|\odr\|_{L^2}^{\frac{2r-p(2-r)}{2p(r-1)}}
\|\nabla\odr\|_{L^2}^{r'\cdot\frac{p-2}{2p}},$$ from which,
$\dH^{\frac 12+\frac 2p}(\R^3)\hookrightarrow
L^{\frac{3p}{p-2}}(\R^3)$ and \eqref{2.3} of Lemma \ref{lem2.4}, we
infer
$$|I_1|\lesssim\int_0^t(\|\u3\|_{\dH^{\frac 12+\frac 2p}}+\|\b3\|_{\dH^{\frac 12+\frac 2p}})\|\odr\|_{L^2}^{\frac 2p}\|\nabla\odr\|_{L^2}^{\frac{2}{p'}}dt'.$$
Applying Young's inequality, we obtain
\begin{align}\label{3.5}
|I_1|\leqslant\frac{r-1}{r^2}\int_0^t\|\nabla\odr\|_{L^2}^2 dt'+
C\int_0^t(\|\u3\|_{\dH^{\frac 12+\frac 2p}}^p+\|\b3\|_{\dH^{\frac
12+\frac 2p}}^p)\|\odr\|_{L^2}^2 dt'.
\end{align}

In order to deal with  $I_2$ and $I_3,$ we need the following lemma:

\begin{lem}[Lemma 4.1 of \cite{CZZ}]\label{lem3.1}
{\sl Let $\theta\in ]0,\alpha(r)[,\sigma\in ]{r'}/{4},1[$, and
$s=\frac 12+1-\frac{2\sigma}{r'}$. Then
\begin{equation}\label{3.6}
\bigl|\int_{\R^3}\ph\Lh f\cdot\ph g\cdot h_{r-1}
dx\bigr|\lesssim\min\bigl\{\|f\|_{L^r},\|f\|_{\htr}\bigr\}\|g\|_{\dH^s}\|h_{\frac
r2}\|_{\dH^\sigma}^{\frac{2}{r'}}.
\end{equation}}
\end{lem}
Next, we estimate $I_2$. We first write by \eqref{1.6}
\begin{equation}\label{3.7}
I_2=I_{2,1}+I_{2,2},
\end{equation}
where
\begin{align*}
& I_{2,1}\eqdefa-\int_0^t\int\bigl(\pa_1{\V_-}\p3\pa_1\Lh{\Gamma_+}+\pa_2{\V_-}\p3\pa_2\Lh{\Gamma_+}\bigr){(\Gamma_+)}_{r-1} dxdt',\\
& I_{2,2}\eqdefa\int_0^t\int\bigl(\pa_1{\V_-}\pa_2\Lh\p3^2{\V_+}-\pa_2{\V_-}\pa_1\Lh\p3^2{\V_+}\bigr){(\Gamma_+)}_{r-1} dxdt'.\\
\end{align*}
Applying Lemma \ref{lem3.1} with
$f=\p3{\Gamma_+},g={\V_-},h=\Gamma_+$,  Gagliardo-Nirenberg
inequality and \eqref{2.3}, we get
\begin{align*}
|I_{2,1}|& \lesssim \int_0^t\|\p3{\Gamma_+}\|_{L^r}\|{\V_-}\|_{\dH^{\frac 32-\frac{2\sigma}{r'}}}\|\odr\|_{\dH^{\sigma}}^{\frac{2}{r'}} dt'\\
& \lesssim \int_0^t\|\nabla\odr\|_{L^2}\|\odr\|_{L^2}^{\frac
2r-1}\|{\V_-}\|_{\dH^{\frac 32-\frac{2\sigma}{r'}}}
\|\odr\|_{L^2}^{\frac{2}{r'}(1-\sigma)}\|\nabla\odr\|_{L^2}^{\frac{2\sigma}{r'}} dt'\\
& \lesssim \int_0^t\|{\V_-}\|_{\dH^{\frac 32-\frac{2\sigma}{r'}}}
\|\odr\|_{L^2}^{2(\frac
12-\frac{\sigma}{r'})}\|\nabla\odr\|_{L^2}^{2(\frac
12+\frac{\sigma}{r'})} dt'.
\end{align*}
Choosing $\sigma=\frac{(p-2)r'}{2p}$, which is between
$\frac{r'}{4}$ and 1 since $p\in \bigl]4,\frac{2r}{2-r}\bigr[$,
gives
$$|I_{2,1}|\lesssim\int_0^t(\|\u3\|_{\dH^{\frac 12+\frac 2p}}+\|\b3\|_{\dH^{\frac 12+\frac 2p}})
\|\odr\|_{L^2}^{\frac 2p}\|\nabla\odr\|_{L^2}^{2(1-\frac 1p)} dt'.$$
Then by using Young's inequality, we get
\begin{equation}\label{3.8}
|I_{2,1}|\leqslant\frac{r-1}{4r^2}\int_0^t\|\nabla\odr\|_{L^2}^2+C\int_0^t(\|\u3\|_{\dH^{\frac
12+\frac 2p}}^p+\|\b3\|_{\dH^{\frac 12+\frac 2p}}^p)
\|\odr\|_{L^2}^2dt'.
\end{equation}
Similarly by applying Lemma \ref{lem3.1} with
$f=\p3^2{\V_+},g={\V_-},h=\Gamma_+$, and
$\sigma=\frac{(p-2)r'}{2p}$, we get
\begin{align*}
|I_{2,2}|& \lesssim \int_0^t\|\p3^2
{\V_+}\|_{\htr}\bigl(\|\u3\|_{\dH^{\frac 32-\frac{2\sigma}{r'}}}
+\|\b3\|_{\dH^{\frac 32-\frac{2\sigma}{r'}}}\bigr)\|\odr\|_{\dH^{\sigma}}^{\frac{2}{r'}}dt'\\
& \lesssim \int_0^t\|\p3^2 {\V_+}\|_{\htr}(\|\u3\|_{\dH^{\frac
12+\frac 2p}}
+\|\b3\|_{\dH^{\frac 12+\frac 2p}})\|\odr\|_{L^2}^{\frac{2}{r'}(1-\sigma)}\|\nabla\odr\|_{L^2}^{\frac{2\sigma}{r'}}dt'\\
& \lesssim \int_0^t\Bigl(\|\p3^2 {\V_+}\|_{\htr}^2\Bigr)^{\frac 12}
\Bigl(\|\u3\|_{\dH^{\frac 12+\frac 2p}}^p+\|\b3\|_{\dH^{\frac 32-\frac{2\sigma}{r'}}}^p\Bigr)^{\alpha(r)}\\
& \qquad\quad\times\Bigl((\|\u3\|_{\dH^{\frac 12+\frac
2p}}^p+\|\b3\|_{\dH^{\frac 12+\frac 2p}}^p)\|\odr\|_{L^2}^2\Bigr)^{\frac
1p-\alpha(r)} \|\nabla\odr\|_{L^2}^{2(\frac 12-\frac 1p)}dt'.
\end{align*}
As we have $\frac 12+\alpha(r)+(\frac 1p-\alpha(r))+(\frac 12-\frac
1p)=1$, applying H\"{o}lder's inequality ensures that
\begin{align*}
|I_{2,2}|\lesssim& \Bigl(\int_0^t \|\p3^2 {\V_+}\|_{\htr}^2
dt'\Bigr)^{\frac 12}
\Bigl(\int_0^t\bigl(\|\u3\|_{\dH^{\frac 12+\frac 2p}}^p+\|\b3\|_{\dH^{\frac 12+\frac 2p}}^p\bigr) dt'\Bigr)^{\alpha(r)}\\
& \times\Bigl(\int_0^t\bigl(\|\u3\|_{\dH^{\frac 12+\frac
2p}}^p+\|\b3\|_{\dH^{\frac 12+\frac 2p}}^p\bigr)\|\odr\|_{L^2}^2
dt'\Bigr)^{\frac 1p-\alpha(r)} \Bigl(\int_0^t\|\nabla\odr\|_{L^2}^2
dt'\Bigr)^{\frac 12-\frac 1p}.
\end{align*}
Then applying Young's inequality leads to
\beq\label{3.9}\begin{split} |I_{2,2}|
\leqslant&\frac{r-1}{4r^2}\int_0^t\|\nabla\odr\|_{L^2}^2 dt'
+C\int_0^t\bigl(\|\u3\|_{\dH^{\frac 12+\frac 2p}}^p+\|\b3\|_{\dH^{\frac 12+\frac 2p}}^p\bigr)\|\odr\|_{L^2}^2 dt'\\
&\qquad\qquad\quad +C\Bigl(\int_0^t\bigl(\|\u3\|_{\dH^{\frac
12+\frac 2p}}^p+\|\b3\|_{\dH^{\frac 12+\frac 2p}}^p\bigr)
dt'\Bigr)^{1-\frac r2} \bigl(\int_0^t \|\p3^2
{\V_+}\|_{\htr}^2dt'\bigr)^{\frac r2}.
\end{split} \eeq
Combining \eqref{3.7}-\eqref{3.9}, we obtain
\begin{align}\label{3.10}
\begin{split}
|I_2|\leqslant& \frac{r-1}{2r^2}\int_0^t\|\nabla\odr\|_{L^2}^2+\|\nabla\omdr\|_{L^2}^2dt'\\
& +C\int_0^t\bigl(\|\u3\|_{\dH^{\frac 12+\frac 2p}}^p+\|\b3\|_{\dH^{\frac 12+\frac 2p}}^p\bigr)\bigl(\|\odr\|_{L^2}^2+\|\omdr\|_{L^2}^2\bigr) dt'\\
& +C\Bigl(\int_0^t\|\u3\|_{\dH^{\frac 12+\frac
2p}}^p+\|\b3\|_{\dH^{\frac 12+\frac 2p}}^p dt'\Bigr)^{1-\frac
r2}\bigl(\int_0^t \|\p3^2 {\V_+}\|_{\htr}^2 dt'\bigr)^{\frac r2}.
\end{split}
\end{align}

We now turn to the last term $I_3$. Use \eqref{1.6} once again, we
write \beq\label{3.11}\begin{split}
I_3 =&2\int_0^t\int\Bigl(\pa_1\bh\cdot\pa_2(\nh^{\bot}\Lh\omega-\nh\Lh\p3\u3)\\
&\qquad\qquad-\pa_2\bh\cdot\pa_1(\nh^{\bot}\Lh\omega-\nh\Lh\p3\u3\Bigr)
{(\Gamma_+)}_{r-1}dxdt'.
\end{split} \eeq
By virtue of  Lemma \ref{lem3.1}, with $f=\nh\omega$,
$g=\bh,\ h=\Gamma_+$, and $\sigma=\frac{(p-2)r'}{2p}$, we get
\begin{align*}
\Bigl|\int_0^t\int\pa_1&\bh\cdot\pa_2\nh^{\bot}\Lh\omega
{(\Gamma_+)}_{r-1}dxdt'\Bigr| \lesssim \int_0^t\|\nh\omega\|_{L^r}\|\bh\|_{\dH^{\frac 12+2(\frac 12-\frac{\sigma}{r'})}}\|\odr\|_{\dH^{\sigma}}^{\frac{2}{r'}} dt'\\
& \lesssim \int_0^t(\|\nh\Gamma_+\|_{L^r}+\|\nh\Gamma_-\|_{L^r})\|\bh\|_{\dH^{\frac 12+2(\frac 12-\frac{\sigma}{r'})}}\|\odr\|_{\dH^{\sigma}}^{\frac{2}{r'}} dt'\\
& \lesssim \int_0^t(\|\nh\odr\|_{L^2}\|\odr\|_{L^2}^{\frac 2r-1}+\|\nh\omdr\|_{L^2}\|\omdr\|_{L^2}^{\frac 2r-1})\\
& \qquad\qquad\times\|\bh\|_{\dH^{\frac 12+2(\frac 12-\frac{\sigma}{r'})}}\|\odr\|_{L^2}^{\frac{2}{r'}(1-\sigma)}\|\nabla\odr\|_{L^2}^{\frac{2\sigma}{r'}} dt'\\
& \lesssim \int_0^t\|\bh\|_{\dH^{\frac 12+\frac
2p}}\bigl(\|\odr\|_{L^2}^{\frac 2p}+\|\omdr\|_{L^2}^{\frac
2p}\bigr)\bigl(\|\nabla\odr\|_{L^2}+\|\nabla\omdr\|_{L^2}\bigr)^{2(1-\frac 1p)} dt',
\end{align*}
where we have used the fact that
$\omega=\frac12\bigl(\Gamma_++\Gamma_-\bigr).$ The same estimate
holds for $
 \int_0^t\int\pa_2\bh\cdot\pa_1\nh^{\bot}\Lh\omega
{(\Gamma_+)}_{r-1}dxdt'.$\\
Along the same line, applying  Lemma \ref{lem3.1} with $f=\nh\p3
u^3,\ g=\bh,\ h=\Gamma_+$, $\sigma=\frac{(p-2)r'}{2p}$,
and the fact that $u^3=\frac12\bigl(\V_++\V_-\bigr)$, yield
\begin{align*}
\Bigl|\int_0^t\int\pa_1\bh\cdot&\pa_2\nh\Lh\p3\u3
{(\Gamma_+)}_{r-1}dxdt'\Bigr| \lesssim \int_0^t\bigl(\|\nh\pa_3u^3\|_{\htr}\|\bh\|_{\dH^{\frac 32-\frac{2\sigma}{r'}}}\|\odr\|_{\dH^{\sigma}}^{\frac{2}{r'}}dt'\\
& \lesssim \int_0^t\|\nh\p3u^3\|_{\htr}\|\bh\|_{\dH^{\frac
12+\frac 2p}}\|\odr\|_{L^2}^{2(\frac 1p-\alpha(r))}
\|\nabla\odr\|_{L^2}^{2(\frac 12-\frac 1p)}dt'\\
&
\lesssim\Bigl(\int_0^t\bigl(\|\nh\p3{\V_+}\|_{\htr}^2+\|\nh\p3{\V_-}\|_{\htr}^2\bigr)dt'\Bigr)^{\frac
12}\Bigl(\int_0^t\|\bh\|_{\dH^{\frac 12+\frac
2p}}^pdt'\Bigr)^{\alpha(r)}\\
&\qquad\qquad\times \Bigl(\int_0^t\|\bh\|_{\dH^{\frac 12+\frac
2p}}^p\|\odr\|_{L^2}^2 dt'\Bigr)^{\frac
1p-\alpha(r)}\Bigl(\int_0^t\|\nabla\odr\|_{L^2}^2dt'\Bigr)^{\frac
12-\frac 1p},
\end{align*}The same estimate holds for
$\int_0^t\int\pa_2\bh\cdot\pa_1\nh\Lh\p3\u3
{(\Gamma_+)}_{r-1}dxdt'.$

 Therefore by applying Young's inequality, we obtain
\begin{align}\label{3.14}
\begin{split}
|I_3|\leqslant&
\frac{r-1}{2r^2}\int_0^t\bigl(\|\nabla\odr\|_{L^2}^2+\|\nabla\omdr\|_{L^2}^2\bigr)
dt'\\
&+C\int_0^t\|\bh\|_{\dH^{\frac 12+\frac 2p}}^p
\bigl(\|\odr\|_{L^2}^2+\|\omdr\|_{L^2}^2\bigr)dt'\\
&+C\Bigl(\int_0^t\|\bh\|_{\dH^{\frac 12+\frac 2p}}^p
dt'\Bigr)^{1-\frac
r2}\Bigl(\int_0^t\bigl(\|\nh\p3{\V_+}\|_{\htr}^2+\|\nh\p3{\V_-}\|_{\htr}^2\bigr)dt'\Bigr)^{\frac
r2}.
\end{split}
\end{align}
Summing up \eqref{3.3}-\eqref{3.5}, \eqref{3.10} and \eqref{3.14}
leads to \beq \label{3.15}
\begin{split}
& \frac 1r\|\odr(t)\|_{L^2}^2+\frac{2(r-1)}{r^2}\int_0^t\|\nabla\odr\|_{L^2}^2-\frac{r-1}{2r^2}\int_0^t\|\nabla\omdr\|_{L^2}^2dt'\\
&\leqslant\frac 1r\bigl(\|\omega_0
\|_{L^r}^r+\|d_0\|_{L^r}^r\bigr)+C\int_0^t\bigl(\|\odr\|_{L^2}^2
+\|\omdr\|_{L^2}^2\bigr)\bigl(\|\u3\|_{\dH^{\frac 12+\frac 2p}}^p+\|b\|_{\dH^{\frac 12+\frac 2p}}^p\bigr)dt'\\
&\ +C\Bigl(\int_0^t
\bigl(\|\nabla\p3{\V_+}\|_{\htr}^2+\|\nabla\p3{\V_-}\|_{\htr}^2\bigr)
dt'\Bigr)^{\frac r2}\Bigl(\int_0^t\bigl(\|\u3\|_{\dH^{\frac 12+\frac
2p}}^p+\|b\|_{\dH^{\frac 12+\frac 2p}}^p\bigr) dt'\Bigr)^{1-\frac
r2}.
\end{split} \eeq
Along the same line, we can get a similar estimate for $\omdr$,
namely \beq \label{3.16}
\begin{split}
& \frac 1r\|\omdr(t)\|_{L^2}^2+\frac{2(r-1)}{r^2}\int_0^t\|\nabla\omdr\|_{L^2}^2-\frac{r-1}{2r^2}\int_0^t\|\nabla\odr\|_{L^2}^2dt'\\
&\leqslant\frac 1r(\|\omega_0 \|_{L^r}^r+\|d_0\|_{L^r}^r)+C\int_0^t\bigl(\|\odr\|_{L^2}^2+\|\omdr\|_{L^2}^2\bigr)
\bigl(\|\u3\|_{\dH^{\frac 12+\frac 2p}}^p+\|b\|_{\dH^{\frac 12+\frac 2p}}^p\bigr)dt'\\
&\ +C\Bigl(\int_0^t
\bigl(\|\nabla\p3{\V_+}\|_{\htr}^2+\|\nabla\p3{\V_-}\|_{\htr}^2\bigr)
dt'\Bigr)^{\frac r2}\Bigl(\int_0^t\bigl(\|\u3\|_{\dH^{\frac 12+\frac
2p}}^p+\|b\|_{\dH^{\frac 12+\frac 2p}}^p\bigr) dt'\Bigr)^{1-\frac
r2}.
\end{split} \eeq
Summing up \eqref{3.15} and \eqref{3.16} and then using  Gronwall's
inequality gives rise to \beq \label{3.18}
\begin{split}
\frac 1r\bigl(&\|\odr(t)
\|_{L^2}^2+\|\omdr(t)\|_{L^2}^2\bigr)+\frac{r-1}{r^2}\int_0^t\bigl(
\|\nabla\odr\|_{L^2}^2+\|\nabla\omdr\|_{L^2}^2\bigr) dt'\\
&\lesssim\Bigl(\frac 2r(\|\omega_0\|_{L^r}^r+\|d_0\|_{L^r}^r)+\bigl(\int_0^t\|\nabla\p3{\V_+}\|_{\htr}^2+\|\nabla\p3{\V_-}\|_{\htr}^2 dt'\bigr)^{\frac r2}\\
&\quad \times\bigl(\int_0^t\|\u3\|_{\dH^{\frac 12+\frac
2p}}^p+\|b\|_{\dH^{\frac 12+\frac 2p}}^p dt'\bigr)^{1-\frac
r2}\Bigr)\cdot\exp\Bigl(C\int_0^t\bigl(\|\u3\|_{\dH^{\frac 12+\frac
2p}}^p+\|b\|_{\dH^{\frac 12+\frac 2p}}^p\bigr) dt'\Bigr).
\end{split} \eeq
This completes the proof of Proposition \ref{prop1.1}, once we
notice the elementary inequality
$$x^{\gamma}e^{c_1x}\lesssim e^{c_2x},\quad \forall \gamma\geqslant0,\ x\geqslant0.$$

\section{Proof of Proposition \ref{prop1.2}}\label{sec4}
Applying $\p3$ on the third components of \eqref{MHD}, we obtain
\beno
\begin{split}
\pt\p3\u3-\Delta\p3\u3=& -\p3 u\cdot\nabla\u3-(u\cdot\nabla)\p3\u3+\p3 b\cdot\nabla \b3+(b\cdot\nabla)\p3\b3\notag\\
&-\p3^2(-\Delta)^{-1}\sum_{\ell,m=1}^3(\pa_\ell u^m\pa_m u^\ell-\pa_\ell b^m\pa_m b^\ell),\\
\pt\p3\b3-\Delta\p3\b3=& -\p3
u\cdot\nabla\b3-(u\cdot\nabla)\p3\b3+\p3 b\cdot\nabla
\u3+(b\cdot\nabla)\p3\u3.
\end{split}
\eeno Adding these two equations gives \beq \label{4.2}
\begin{split}
\pt\p3\V_+-\Delta\p3\V_+=&-\p3 u\cdot\nabla{\V_+}+\p3 b\cdot\nabla{\V_+}-u\cdot\nabla\p3\V_+\\
&+b\cdot\nabla\p3\V_+-\p3^2(-\Delta)^{-1}\sum_{\ell,m=1}^3\bigl(\pa_\ell
u^m\pa_m u^\ell-\pa_\ell b^m\pa_m b^\ell\bigr).
\end{split} \eeq
We write
\begin{align*}
-\p3 u\cdot\nabla{\V_+}+\p3 b\cdot\nabla {\V_+}=-\sum_{\ell=1}^2
\p3(u^\ell-b^\ell)\pa_\ell{\V_+}-(\p3\u3)^2+(\p3\b3)^2,
\end{align*}
and
\begin{align*}
 \sum_{\ell,m=1}^3(\pa_\ell u^m\pa_m u^\ell-\pa_\ell b^m\pa_m b^\ell)
=& \sum_{\ell,m=1}^2(\pa_\ell u^m\pa_m u^\ell-\pa_\ell b^m\pa_m
b^\ell)\\
&+2\sum_{\ell=1}^2(\p3 u^\ell\pa_\ell\u3-\p3
b^\ell\pa_\ell\b3)+(\p3\u3)^2-(\p3\b3)^2.
\end{align*}
Then we take the $\htr$ inner product of \eqref{4.2} with
$\p3{\V_+}$ to obtain
\begin{align}\label{4.3}
\frac 12\frac{d}{dt}\|\p3&
{\V_+}(t)\|_{\htr}^2+\|\nabla\p3{\V_+}\|_{\htr}^2=-\sum_{i=1}^5\Rmnum{2}_i,
\end{align}
with
\begin{align*}
& \Rmnum{2}_1\eqdefa\bigl((Id+\p3^2(-\Delta)^{-1})((\p3\u3)^2-(\p3\b3)^2)\big|\p3{\V_+}\bigr)_{\htr},\\
& \Rmnum{2}_2\eqdefa\bigl(\p3^2(-\Delta)^{-1}\sum_{\ell,m=1}^2\pa_\ell u^m\pa_m u^\ell-\pa_\ell b^m\pa_m b^\ell\big|\p3{\V_+}\bigr)_{\htr},\\
& \Rmnum{2}_3\eqdefa\sum_{\ell=1}^2\bigl(\p3(u^\ell-b^\ell)\pa_\ell{\V_+}\big|\p3{\V_+}\bigr)_{\htr},\\
& \Rmnum{2}_4\eqdefa\bigl(2\p3^2(-\Delta)^{-1}\sum_{\ell=1}^2(\p3 u^\ell\pa_\ell\u3-\p3 b^\ell\pa_\ell\b3)\big|\p3{\V_+}\bigr)_{\htr},\\
&
\Rmnum{2}_5\eqdefa\bigl((u-b)\cdot\nabla\p3{\V_+}\big|\p3{\V_+}\bigr)_{\htr}.
\end{align*}
Let us first recall the following lemma from \cite{CZZ}.
\begin{lem}[Lemma 5.1 of \cite{CZZ}]\label{lem4.1}
{\ Let A be a bounded Fourier multiplier. If p
and $\theta$ satisfy
\begin{equation}\label{4.4}
0<\theta<\frac 12-\frac 1p,
\end{equation}
then we have
$$\bigl|\bigl(A(D)(fg)\big|\p3 h^3\bigr)_{\htr}\bigr|\lesssim\|f\|_{\dH^{\frac 12-3\alpha(r)+\theta,\frac 12-\frac 1p-\theta}}
\|g\|_{\dH^{\frac 12-3\alpha(r)+\theta,\frac 12-\frac
1p-\theta}}\|h^3\|_{\dH^{\frac 12+\frac 2p}}.$$}
\end{lem}
Noting that $p>4,r>\frac 43$, we have $\frac 1p+\frac 1r<1$, and
hence $\theta<\alpha(r)<\frac 12-\frac 1p$, i.e.  the condition
\eqref{4.4} is satisfied under the assumption of Proposition
\ref{prop1.2}. Because $\p3^2\Delta^{-1}$ is a bounded Fourier
multiplier, applying Lemma \ref{lem4.1} with
$f=\p3\V_+,g=\p3\V_-,h=u+b$, gives
\begin{align}\label{4.5}
|\Rmnum{2}_1|\lesssim& \|\V_+\|_{\dH^{\frac 12+\frac
2p}}\|\p3{\V_+}\|_{\dH^{\frac 12-3\alpha(r)+\theta,\frac 12-\frac
1p-\theta}}\|\p3{\V_-}\|_{\dH^{\frac 12-3\alpha(r)+\theta,\frac
12-\frac 1p-\theta}}.
\end{align}
 While a direct calculation from Definition
\ref{def1.2} gives
\begin{align*}
\|a\|_{\dH^{\frac 12-3\alpha(r)+\theta,\frac 12-\frac 1p-\theta}}^2
& \leqslant\int_{\R^3}|\widehat{a}(\xi)|^{\frac 2p}\bigl(|\xi||\widehat{a}(\xi)|\bigr)^{\frac{2}{p'}}|\xih|^{2(-3\alpha(r)+\theta)}|\xi_3|^{-2\theta}d\xi\\
& \leqslant\Bigl(\int_{\R^3}|\widehat{a}(\xi)|^2|\xih|^{2(-3\alpha(r)+\theta)}|\xi_3|^{-2\theta}d\xi\Bigr)^{\frac 1p}\\
&
\qquad\times\Bigl(\int_{\R^3}|\widehat{a}(\xi)|^2|\xi|^2|\xih|^{2(-3\alpha(r)+\theta)}|\xi_3|^{-2\theta}d\xi\Bigr)^{\frac{1}{p'}},
\end{align*}
that is
\begin{equation}\label{4.7}
\|a\|_{\dH^{\frac 12-3\alpha(r)+\theta,\frac 12-\frac
1p-\theta}}\leqslant\|a\|_{\htr}^{\frac 1p}\|\nabla
a\|_{\htr}^{\frac {1}{p'}}.
\end{equation}
Applying \eqref{4.7} to $\p3\V_\pm$ in \eqref{4.5} gives
\begin{align*}
|\Rmnum{2}_1|\lesssim\|\V_+\|_{\dH^{\frac 12+\frac 2p}}
\|\p3{\V_+}\|_{\htr}^{\frac 1p}\|\nabla \p3{\V_+}\|_{\htr}^{\frac
{1}{p'}}\|\p3{\V_-}\|_{\htr}^{\frac 1p}\|\nabla
\p3{\V_-}\|_{\htr}^{\frac {1}{p'}},
\end{align*}
then Young's inequality and mean inequality ensure
\begin{align}\label{4.8}
\begin{split}
|\Rmnum{2}_1|\leqslant& \frac {1}{20}\bigl(\|\nabla
\p3{\V_+}\|_{\htr}^2+\|\nabla
\p3{\V_-}\|_{\htr}^2\bigr)+C\|\V_+\|_{\dH^{\frac 12+\frac
2p}}^p\bigl(\|\p3{\V_+}\|_{\htr}^2+\|\p3{\V_-}\|_{\htr}^2\bigr).
\end{split}
\end{align}

For $\Rmnum{2}_2$ term, one key point is to write it in a symmetric
form, namely
$$\Rmnum{2}_2=\Bigl(\p3^2(-\Delta)^{-1}\sum_{\ell,m=1}^2\frac 12 \bigl(\pa_\ell (u^m+b^m)\pa_m (u^\ell-b^\ell)+\pa_\ell (u^m-b^m)\pa_m (u^\ell+b^\ell)\bigr)|\p3{\V_+}\Bigr)_{\htr}.$$
Applying the Hodge decomposition for the horizontal variables to
$u^\ell\pm b^\ell$, and noting both $\p3^2\Delta^{-1}$ and
$\ph^2\Delta_\h^{-1}$ are bounded Fourier multipliers, then Lemma
\ref{lem4.1} ensures that
\begin{align*}
|\Rmnum{2}_2|\lesssim&\|\V_+\|_{\dH^{\frac 12+\frac 2p}}
\bigl(\|\p3{\V_+}\|_{\dH^{\frac 12-3\alpha(r)+\theta,\frac 12-\frac 1p-\theta}}+\|\Gamma_+\|_{\dH^{\frac 12-3\alpha(r)+\theta,\frac 12-\frac 1p-\theta}}\bigr)\\
&\qquad\qquad\times\bigl(\|\p3{\V_-}\|_{\dH^{\frac
12-3\alpha(r)+\theta,\frac 12-\frac
1p-\theta}}+\|\Gamma_-\|_{\dH^{\frac 12-3\alpha(r)+\theta,\frac
12-\frac 1p-\theta}}\bigr).
\end{align*}
Yet it follows from  Lemma \ref{lem2.1} and Lemma \ref{lem2.4} that
for any function a
\begin{equation}\label{4.6}
\|a\|_{\dH^{\frac 12-3\alpha(r)+\theta,\frac 12-\frac
1p-\theta}}\leqslant\|a\|_{\dH^{1-3\alpha(r)-\frac 1p}}
\lesssim\|a_{\frac r2}\|_{L^2}^{2\alpha(r)+\frac 1p}\|\nabla
a_{\frac r2}\|_{L^2}^{\frac {1}{p'}},
\end{equation}
where $p'$ denotes the conjugate index of p. Then applying \eqref{4.6} to $\Gamma_\pm$,
\eqref{4.7} to $\p3\V_\pm$ gives
\begin{align*}
|\Rmnum{2}_2|\lesssim&\|\V_+\|_{\dH^{\frac 12+\frac 2p}}
\bigl(\|\p3{\V_+}\|_{\htr}^{\frac 1p}\|\nabla \p3{\V_+}\|_{\htr}^{\frac {1}{p'}}+\|\odr\|_{L^2}^{\frac 1p+2\alpha(r)}\|\nabla \odr\|_{L^2}^{\frac {1}{p'}}\bigr)\\
&\times\bigl(\|\p3{\V_-}\|_{\htr}^{\frac 1p}\|\nabla
\p3{\V_-}\|_{\htr}^{\frac {1}{p'}}+\|\omdr\|_{L^2}^{\frac
1p+2\alpha(r)}\|\nabla \omdr\|_{L^2}^{\frac {1}{p'}}\bigr),
\end{align*}
which implies that \beq\label{4.9}\begin{split}
|&\Rmnum{2}_2|\leqslant \frac {1}{20}\bigl(\|\nabla
\p3{\V_+}\|_{\htr}^2+\|\nabla \p3{\V_-}\|_{\htr}^2\bigr)
+C\|\V_+\|_{\dH^{\frac 12+\frac 2p}}^p\bigl(\|\p3{\V_+}\|_{\htr}^2+\|\p3{\V_-}\|_{\htr}^2\bigr)\\
&\ \ +C\|\V_+\|_{\dH^{\frac 12+\frac
2p}}\Bigl(\|\odr\|_{L^2}^{2(\frac 1p+2\alpha(r))}\|\nabla
\odr\|_{L^2}^{\frac {2}{p'}}+\|\omdr\|_{L^2}^{2(\frac
1p+2\alpha(r))}\|\nabla \omdr\|_{L^2}^{\frac {2}{p'}}\Bigr).
\end{split} \eeq


On the other hand, for any real valued functions a and b, and any
couple $(\alpha,\beta)\in\R^2$, applying H\"{o}lder's inequality
gives
\begin{align}\label{4.10}
\begin{split}
|(a|b)_{\htr}|&
=\Bigl|\int_{\R^3}|\xih|^{-6\alpha(r)+2\theta-\alpha}|\xi_3|^{-\beta-2\theta}\widehat{a}(\xi)
|\xih|^{\alpha}|\xi_3|^{\beta}\widehat{b}(\xi)d\xi\Bigr|\\
&\leqslant\|a\|_{\dH^{-6\alpha(r)+2\theta-\alpha,-\beta-2\theta}}\|b\|_{\dH^{\alpha,\beta}}.
\end{split}
\end{align}
Note that $4<p<\frac{2r}{2-r}$ and $3\alpha(r)-\frac
2p<\theta<\alpha(r)$, we have $\frac
2p+3\alpha(r)-\theta\in]0,1[$, and hence
\begin{align*}
|\xih|^{2(1-6\alpha(r)-\frac 2p+2\theta)}|\xi_3|^{2(\frac 2p+3\alpha(r)-2\theta)}& =|\xih|^{2(-3\alpha(r)+\theta)}|\xi_3|^{-2\theta}\cdot|\xih|^{2(1-3\alpha(r)-\frac 2p+\theta)}|\xi_3|^{2(\frac 2p+3\alpha(r)-\theta)}\\
&\leqslant|\xih|^{2(-3\alpha(r)+\theta)}|\xi_3|^{-2\theta}|\xi|^2,
\end{align*}
which implies for any function a
\begin{align}\label{4.12}
\begin{split}
\|a\|_{\dH^{1-6\alpha(r)-\frac 2p+2\theta,\frac 2p+3\alpha(r)-2\theta}}^2& =\int_{\R^3}|\xih|^{2(1-6\alpha(r)-\frac 2p+2\theta)}|\xi_3|^{2(\frac 2p+3\alpha(r)-2\theta)}|\widehat{a}(\xi)|^2d\xi\\
& \leqslant\int_{\R^3}|\xih|^{2(-3\alpha(r)+\theta)}|\xi_3|^{-2\theta}\bigl(|\xi||\widehat{a}(\xi)|\bigr)^2d\xi\\
& =\|\nabla a\|_{\htr}^2.
\end{split}
\end{align}
Along the same line, one has
\begin{align}\label{4.13}
\|a\|_{\dH^{1-6\alpha(r)+2\theta,3\alpha(r)-2\theta}}\leqslant\|\nabla
a\|_{\htr}.
\end{align}
In order to estimate $\Rmnum{2}_3$, we apply Bony's decomposition in
the vertical variable to write
$$\p3 (u^\ell-b^\ell)\pa_\ell \V_+=(T^\v+\bar{T}^\v+R^\v)\bigl(\p3 (u^\ell-b^\ell),\pa_\ell \V_+\bigr).$$
Applying \eqref{4.10} with $\alpha=1-6\alpha(r)-\frac
2p+2\theta,\beta=\frac 2p+3\alpha(r)-2\theta$, the law of
product (see Lemma 4.5 of \cite{CZ14} for
example) and \eqref{4.12} ensures that
\begin{align}\label{4.15}
\begin{split}
& \Bigl|\Bigl(\sum_{\ell=1}^2 (T^\v+\bar{T}^\v)\bigl(\p3 (u^\ell-b^\ell),\pa_{\ell}\V_+\bigr)\Big|\p3\V_+\Bigr)_{\htr}\Bigr|\\
& \lesssim\sum_{\ell=1}^2 \|(T^\v+\bar{T}^\v)\bigl(\p3 (u^\ell-b^\ell),\pa_\ell \V_+\bigr)\|_{\dH^{\frac 2p-1,-3\alpha(r)-\frac 2p}}\|\p3\V_+\|_{\dH^{1-6\alpha(r)-\frac 2p+2\theta,\frac 2p+3\alpha(r)-2\theta}}\\
& \lesssim\sum_{\ell=1}^2 \|\p3 (u^\ell-b^\ell)\|_{(\dB^1_{2,1})_\h(\dB^{-3\alpha(r)-\frac
2p}_{2,1})_\v}\|\pa_\ell \V_+\|_{(\dH^{\frac 2p-1})_\h(\dB^{\frac
12}_{2,1})_\v}\|\nabla\p3 {\V_+}\|_{\htr}\\
&\lesssim\sum_{\ell=1}^2\|u^\ell-b^\ell\|_{(\dB^1_{2,1})_\h(\dB^{1-3\alpha(r)-\frac
2p}_{2,1})_\v}\|\V_+\|_{\dH^{\frac
12+\frac 2p}}\|\nabla\p3 {\V_+}\|_{\htr}.
\end{split}
\end{align}
Applying \eqref{4.10} with $\alpha=0,\ \beta=-\frac 12+\frac
2p$, the law of product, and \eqref{4.13} ensures that
\begin{align}\label{4.16}
\begin{split}
& \Bigl|\Bigl(\sum_{\ell=1}^2 R^\v\bigl(\p3 (u^\ell-b^\ell),\pa_{\ell}\V_+\bigr)\Big|\p3\V_+\Bigr)_{\htr}\Bigr|\\
&\lesssim\sum_{\ell=1}^2\|R^\v\bigl(\p3 (u^\ell-b^\ell),\pa_\ell {\V_+}\bigr)\|_{\dH^{-6\alpha(r)+2\theta,\frac 12-\frac 2p-2\theta}}\|\p3{\V_+}\|_{\dH^{0,-\frac 12+\frac 2p}}\\
&\lesssim\sum_{\ell=1}^2\|\p3 (u^\ell-b^\ell)\|_{(\dB^1_{2,1})_\h(\dB^{-3\alpha(r)-\frac 2p}_{2,1})_\v}
\|\pa_\ell {\V_+}\|_{\dH^{-6\alpha(r)+2\theta,1+3\alpha(r)-2\theta}}\|\V_+\|_{\dH^{0,\frac 12+\frac 2p}}\\
&\lesssim\sum_{\ell=1}^2\|u^\ell-b^\ell\|_{(\dB^1_{2,1})_\h(\dB^{1-3\alpha(r)-\frac
2p}_{2,1})_\v}\|\nabla\p3 {\V_+}\|_{\htr}\|\V_+\|_{\dH^{\frac
12+\frac 2p}}.
\end{split}
\end{align}
Therefore, by virtue of Lemma \ref{lem2.5} with $\beta=\frac 2p$,
inequalities \eqref{4.15} and \eqref{4.16} ensure that
\begin{align*}
|\Rmnum{2}_3|
\lesssim& \|\V_+\|_{\dH^{\frac 12+\frac 2p}}\|\nabla\p3 {\V_+}\|_{\htr}\\
& \times\bigl(\|\omdr\|_{L^2}^{2(\alpha(r)+\frac
1p)}\|\nabla\omdr\|_{L^2}^{1-\frac 2p} +\|\p3{\V_-}\|_{\htr}^{\frac
2p}\|\nabla\p3{\V_-}\|_{\htr}^{1-\frac 2p}\bigr).
\end{align*}
Applying Young's inequality and mean inequality yields
\begin{align}\label{4.17}
\begin{split}
|\Rmnum{2}_3|
\leqslant& \frac {1}{20}\bigl(\|\nabla\p3 {\V_+}\|_{\htr}^2+\|\nabla\p3 {\V_-}\|_{\htr}^2\bigr)+C\|\V_+\|_{\dH^{\frac 12+\frac 2p}}^p\|\p3{\V_-}\|_{\htr}^2\\
& +C\|\V_+\|_{\dH^{\frac 12+\frac
2p}}^2\|\omdr\|_{L^2}^{4(\alpha(r)+\frac
1p)}\|\nabla\omdr\|_{L^2}^{2(1-\frac 2p)}.
\end{split}
\end{align}

The term $\Rmnum{2}_4$ can be handled as above. Indeed we first
rewrite it in a symmetric form
$$\Rmnum{2}_4=\Bigl(2\p3^2(-\Delta)^{-1}\sum_{\ell=1}^2\bigl(\p3 (u^\ell+b^\ell)\pa_\ell{\V_-}+\p3 (u^\ell-b^\ell)\pa_\ell{\V_+}\bigr)\Big|\p3{\V_+}\Bigr)_{\htr}.$$
Then it follows from the estimate of  $\Rmnum{2}_3$ that
\begin{align*}
|\Rmnum{2}_4|
\lesssim& \bigl(\|\V_+\|_{\dH^{\frac 12+\frac 2p}}+\|{\V_-}\|_{\dH^{\frac 12+\frac 2p}}\bigr)\bigl(\|\nabla\p3 {\V_+}\|_{\htr}+\|\nabla\p3 {\V_-}\|_{\htr}\bigr)\\
& \times\Bigl(\|\odr\|_{L^2}^{2(\alpha(r)+\frac
1p)}\|\nabla\odr\|_{L^2}^{1-\frac 2p}
+\|\p3{\V_+}\|_{\htr}^{\frac 2p}\|\nabla\p3{\V_+}\|_{\htr}^{1-\frac 2p}\\
&\quad +\|\omdr\|_{L^2}^{2(\alpha(r)+\frac
1p)}\|\nabla\omdr\|_{L^2}^{1-\frac 2p} +\|\p3{\V_-}\|_{\htr}^{\frac
2p}\|\nabla\p3{\V_-}\|_{\htr}^{1-\frac 2p}\Bigr).
\end{align*}
Applying Young's inequality  yields
\begin{align}\label{4.18}
|\Rmnum{2}_4|
\leqslant\frac {1}{20}& \bigl(\|\nabla\p3 {\V_+}\|_{\htr}^2+\|\nabla\p3 {\V_-}\|_{\htr}^2\bigr)+C\bigl(\|\V_+\|_{\dH^{\frac 12+\frac 2p}}^2+\|{\V_-}\|_{\dH^{\frac 12+\frac 2p}}^2\bigr)\notag\\
\times&\bigl(\|\odr\|_{L^2}^{4(\alpha(r)+\frac 1p)}\|\nabla\odr\|_{L^2}^{2(1-\frac 2p)}+\|\omdr\|_{L^2}^{4(\alpha(r)+\frac 1p)}\|\nabla\omdr\|_{L^2}^{2(1-\frac 2p)}\bigr)\\
+&C\bigl(\|\V_+\|_{\dH^{\frac 12+\frac 2p}}^p+\|{\V_-}\|_{\dH^{\frac
12+\frac
2p}}^p\bigr)\bigl(\|\p3{\V_+}\|_{\htr}^2+\|\p3{\V_-}\|_{\htr}^2\bigr).\notag
\end{align}

Finally let us turn to the estimate of $\Rmnum{2}_5$. We first
decompose it as
\begin{align}\label{4.19}
\begin{split}
\Rmnum{2}_5=&\bigl((\uh-\bh)\cdot\nh\pa_3{\V_+}\big|\p3{\V_+}\bigr)_{\htr}+\bigl({\V_-}\cdot\p3^2{\V_+}\big|\p3{\V_+}\bigr)_{\htr}
\eqdefa\Rmnum{2}_{5,1}+\Rmnum{2}_{5,2}.
\end{split}
\end{align}
Applying Lemma \ref{lem2.6} gives
\begin{align}\label{4.20}
\begin{split}
|\Rmnum{2}_{5,1}|\lesssim &\Bigl(\|\nh(\uh-\bh)\|_{\dH^{\frac
12-3\alpha(r)+\theta,\frac 12-\frac
1p-\theta}}^2+\|\p3{\V_+}\|_{\dH^{\frac 12-3\alpha(r)+\theta,\frac
12-\frac 1p-\theta}}^2
\\
&\quad+\|\uh-\bh\|_{(\dB^1_{2,1})_\h(\dB^{1-3\alpha(r)-\frac
2p}_{2,1})_\v}\|\nabla\p3{\V_+}\|_{\htr}\Bigr)\cdot\|\V_+\|_{\dH^{\frac
12+\frac 2p}}.
\end{split}
\end{align}
While using Hodge decomposition
\eqref{1.6}, and then \eqref{4.7},\eqref{4.6}, we have
\begin{align*}
\|\nh&(\uh-\bh)\|_{\dH^{\frac 12-3\alpha(r)+\theta,\frac 12-\frac
1p-\theta}}^2
+\|\p3{\V_+}\|_{\dH^{\frac 12-3\alpha(r)+\theta,\frac 12-\frac 1p-\theta}}^2\\
\lesssim&\|\Gamma_-\|_{\dH^{\frac 12-3\alpha(r)+\theta,\frac 12-\frac 1p-\theta}}^2+\|\p3{\V_-}\|_{\dH^{\frac 12-3\alpha(r)+\theta,\frac 12
-\frac 1p-\theta}}^2+\|\p3{\V_+}\|_{\dH^{\frac 12-3\alpha(r)+\theta,\frac 12-\frac 1p-\theta}}^2\\
\lesssim& \|\omdr\|_{L^2}^{2(2\alpha(r)+\frac
1p)}\|\nabla\omdr\|_{L^2}^{\frac{2}{p'}}+\|\p3{\V_-}\|_{\htr}^{\frac
2p}\|\nabla\p3{\V_-}\|_{\htr}^{\frac {2}{p'}}
+\|\p3{\V_+}\|_{\htr}^{\frac 2p}\|\nabla\p3{\V_+}\|_{\htr}^{\frac
{2}{p'}}.
\end{align*}
Inserting this estimate and \eqref{2.5} with $\beta=\frac 2p$ into
\eqref{4.20} yields
\begin{align}\label{4.21}
\begin{split}
|\Rmnum{2}_{5,1}|\lesssim& \|\V_+\|_{\dH^{\frac 12+\frac 2p}}\Bigl(\|\omdr\|_{L^2}^{2(2\alpha(r)+\frac 1p)}\|\nabla\omdr\|_{L^2}^{\frac{2}{p'}}
 +\|\p3{\V_+}\|_{\htr}^{\frac 2p}\|\nabla\p3{\V_+}\|_{\htr}^{\frac {2}{p'}}\\
& +\|\p3{\V_-}\|_{\htr}^{\frac 2p}\|\nabla\p3{\V_-}\|_{\htr}^{\frac {2}{p'}}+\|\nabla\p3{\V_+}\|_{\htr}\\
& \times \bigl(\|\omdr\|_{L^2}^{2(\alpha(r)+\frac
1p)}\|\nabla\omdr\|_{L^2}^{1-\frac 2p} +\|\p3{\V_-}\|_{\htr}^{\frac
2p}\|\nabla\p3{\V_-}\|_{\htr}^{1-\frac 2p}\bigr)\Bigr).
\end{split}
\end{align}

In order to estimate $\Rmnum{2}_{5,2}$, we need the following lemma:
\begin{lem}[(95) of \cite{Yamazaki}]\label{lem4.2}
{\sl Let $s_1<1,s_2<1,s_1+s_2>0$ and $0<\sigma<1$. Then we
have
\begin{equation}\label{4.22}
\|fg\|_{\dH^{s_1+s_2-1,\sigma-\frac
12}}\lesssim\|f\|_{(\dB^{s_1}_{2,2})_\h(\dB^{\frac
12}_{2,1})_\v}\|g\|_{\dH^{s_2,\sigma-\frac 12}}.
\end{equation}}
\end{lem}
This lemma together with \eqref{4.10} and Lemma \ref{lem2.1} ensure
that
\begin{align}\label{4.23}
\begin{split}
|\Rmnum{2}_{5,2}|& \leqslant\|{\V_-}\p3^2{\V_+}\|_{\dH^{-1-3\alpha(r)+\frac 2p+\theta,-\theta}}\|\p3{\V_+}\|_{\dH^{1-3\alpha(r)-\frac 2p+\theta,-\theta}}\\
& \lesssim\|{\V_-}\|_{(\dB^{\frac 2p}_{2,2})_\h(\dB^{\frac 12}_{2,1})_\v}\|\p3^2{\V_+}\|_{\htr}\|\p3{\V_+}\|_{\dH^{1-3\alpha(r)-\frac 2p+\theta,-\theta}}\\
& \lesssim\|{\V_-}\|_{\dH^{\frac 12+\frac
2p}}\|\p3^2{\V_+}\|_{\htr}\|\p3{\V_+}\|_{\dH^{1-3\alpha(r)-\frac
2p+\theta,-\theta}}.
\end{split}
\end{align}
This along with the interpolation, which claims that for any
function a
\begin{align}\label{4.24}
\begin{split}
\|a\|_{\dH^{1-3\alpha(r)-\frac 2p+\theta,-\theta}}^2
& \leqslant\Bigl(\int_{\R^3}|\xih|^{-6\alpha(r)+2\theta}|\xi_3|^{-2\theta}|\widehat{a}(\xi)|^2d\xi\Bigr)^{\frac 2p}\\
& \quad\times\Bigl(\int_{\R^3}|\xih|^2|\xih|^{-6\alpha(r)+2\theta}|\xi_3|^{-2\theta}|\widehat{a}(\xi)|^2d\xi\Bigr)^{1-\frac 2p}\\
& =\|a\|_{\htr}^{\frac 4p}\|\nh a\|_{\htr}^{2(1-\frac 2p)},
\end{split}
\end{align}
ensures that
\begin{align}\label{4.25}
\begin{split}
|\Rmnum{2}_{5,2}|& \lesssim \|{\V_-}\|_{\dH^{\frac 12+\frac 2p}}\|\p3^2{\V_+}\|_{\htr}\|\p3{\V_+}\|_{\htr}^{\frac 2p}\|\nh \p3{\V_+}\|_{\htr}^{1-\frac 2p}\\
&\lesssim \|{\V_-}\|_{\dH^{\frac 12+\frac
2p}}\|\p3{\V_+}\|_{\htr}^{\frac 2p}\|\nabla
\p3{\V_+}\|_{\htr}^{\frac {2}{p'}}.
\end{split}
\end{align}
Hence, by summing up  \eqref{4.21} and \eqref{4.25}, we obtain
\begin{align*}
|\Rmnum{2}_5|\lesssim& \Bigl(\|\V_+\|_{\dH^{\frac 12+\frac 2p}}+\|{\V_-}\|_{\dH^{\frac 12+\frac 2p}}\Bigr)\Bigl(\|\omdr\|_{L^2}^{2(2\alpha(r)+\frac 1p)}\|\nabla\omdr\|_{L^2}^{\frac{2}{p'}}\notag\\
& +\|\p3{\V_+}\|_{\htr}^{\frac 2p}\|\nabla\p3{\V_+}\|_{\htr}^{\frac {2}{p'}}+\|\p3{\V_-}\|_{\htr}^{\frac 2p}\|\nabla\p3{\V_-}\|_{\htr}^{\frac {2}{p'}}\notag\\
&
+\|\nabla\p3{\V_+}\|_{\htr}\bigl(\|\omdr\|_{L^2}^{2(\alpha(r)+\frac
1p)}\|\nabla\omdr\|_{L^2}^{1-\frac 2p}\notag
+\|\p3{\V_-}\|_{\htr}^{\frac 2p}\|\nabla\p3{\V_-}\|_{\htr}^{1-\frac
2p}\bigr)\Bigr),
\end{align*} which implies
\begin{align}\label{4.26}
\begin{split}
|\Rmnum{2}_5|\leqslant&\frac {1}{20}\bigl(\|\nabla\p3{\V_+}\|_{\htr}^2+\|\nabla\p3{\V_-}\|_{\htr}^2\bigr)\\
&+C\bigl(\|\u3\|_{\dH^{\frac 12+\frac 2p}}+\|\b3\|_{\dH^{\frac 12+\frac 2p}}\bigr)\bigl(\|\omdr\|_{L^2}^{2(2\alpha(r)+\frac 1p)}\|\nabla\omdr\|_{L^2}^{\frac{2}{p'}}\bigr)\\
&+C\bigl(\|\u3\|_{\dH^{\frac 12+\frac 2p}}^p+\|\b3\|_{\dH^{\frac 12+\frac 2p}}^p\bigr)\bigl(\|\p3{\V_+}\|_{\htr}^2+\|\p3{\V_-}\|_{\htr}^2\bigr)\\
&+C\bigl(\|\u3\|_{\dH^{\frac 12+\frac 2p}}^2+\|\b3\|_{\dH^{\frac 12+\frac
2p}}^2\bigr)\bigl(\|\omdr\|_{L^2}^{4(\alpha(r)+\frac
1p)}\|\nabla\omdr\|_{L^2}^{2(1-\frac 2p)}\bigr).
\end{split}
\end{align}
Substituting the estimates \eqref{4.8}, \eqref{4.9}, \eqref{4.17},
\eqref{4.18} and \eqref{4.26} into \eqref{4.3} leads to
\begin{align}\label{4.27}
\begin{split}
 \frac
12\frac{d}{dt}&\|\p3{\V_+}(t)\|_{\htr}^2+\|\nabla\p3{\V_+}\|_{\htr}^2
\leqslant \frac {1}{4}\bigl(\|\nabla\p3{\V_+}\|_{\htr}^2+\|\nabla\p3{\V_-}\|_{\htr}^2\bigr)\\
& + C\bigl(\|\u3\|_{\dH^{\frac 12+\frac 2p}}^p+\|\b3\|_{\dH^{\frac 12+\frac 2p}}^p\bigr)\bigl(\|\p3{\V_+}\|_{\htr}^2+\|\p3{\V_-}\|_{\htr}^2\bigr)\\
& + C\bigl(\|\u3\|_{\dH^{\frac 12+\frac 2p}}+\|\b3\|_{\dH^{\frac 12+\frac 2p}}\bigr)\Bigl(\|\odr\|_{L^2}^{2(2\alpha(r)+\frac 1p)}\|\nabla\odr\|_{L^2}^{\frac{2}{p'}}\\
& +\|\omdr\|_{L^2}^{2(2\alpha(r)+\frac 1p)}\|\nabla\omdr\|_{L^2}^{\frac{2}{p'}}\Bigr)+C\bigl(\|\u3\|_{\dH^{\frac 12+\frac 2p}}^2+\|\b3\|_{\dH^{\frac 12+\frac 2p}}^2\bigr)\\
\times& \bigl(\|\odr\|_{L^2}^{4(\alpha(r)+\frac
1p)}\|\nabla\odr\|_{L^2}^{2(1-\frac
2p)}+\|\omdr\|_{L^2}^{4(\alpha(r)+\frac
1p)}\|\nabla\omdr\|_{L^2}^{2(1-\frac 2p)}\bigr).
\end{split}
\end{align}
Exactly along the same line to the derivation of the above
inequality, we have
\begin{align}\label{4.28}
\begin{split}
\frac
12\frac{d}{dt}&\|\p3{\V_-}(t)\|_{\htr}^2+\|\nabla\p3{\V_-}\|_{\htr}^2
\leqslant \text{the right hand side of \eqref{4.27}}.
\end{split}
\end{align}
Summing up the above two estimates gives rise to
\begin{align}\label{4.29}
\begin{split}
 \frac{d}{dt}\bigl(&\|\p3{\V_+}(t)\|_{\htr}^2+\|\p3{\V_-}(t)\|_{\htr}^2\bigr)+\|\nabla\p3{\V_+}\|_{\htr}^2+\|\nabla\p3{\V_-}\|_{\htr}^2\\
\leqslant& C\bigl(\|\u3\|_{\dH^{\frac 12+\frac 2p}}^p+\|\b3\|_{\dH^{\frac 12+\frac 2p}}^p\bigr)\bigl(\|\p3{\V_+}\|_{\htr}^2+\|\p3{\V_-}\|_{\htr}^2\bigr)\\
& +C\bigl(\|\u3\|_{\dH^{\frac 12+\frac 2p}}+\|\b3\|_{\dH^{\frac 12+\frac 2p}}\bigr)\Bigl(\|\odr\|_{L^2}^{2(2\alpha(r)+\frac 1p)}\|\nabla\odr\|_{L^2}^{\frac{2}{p'}}\\
& +\|\omdr\|_{L^2}^{2(2\alpha(r)+\frac 1p)}\|\nabla\omdr\|_{L^2}^{\frac{2}{p'}}\Bigr)+C\bigl(\|\u3\|_{\dH^{\frac 12+\frac 2p}}^2+\|\b3\|_{\dH^{\frac 12+\frac 2p}}^2\bigr)\\
\times& \bigl(\|\odr\|_{L^2}^{4(\alpha(r)+\frac
1p)}\|\nabla\odr\|_{L^2}^{2(1-\frac
2p)}+\|\omdr\|_{L^2}^{4(\alpha(r)+\frac
1p)}\|\nabla\omdr\|_{L^2}^{2(1-\frac 2p)}\bigr).
\end{split}
\end{align}
Then Gronwall's inequality allows to conclude the proof of
Proposition \ref{prop1.2} by noticing that
$$\|\p3\u3(0)\|_{\htr}^2+\|\p3\b3(0)\|_{\htr}^2\lesssim\|u(0)\|_{\dH^{1-3\alpha(r)}}^2+\|b(0)\|_{\dH^{1-3\alpha(r)}}^2
\lesssim\|\Omega_0\|_{L^r}^2+\|j_0\|_{L^r}^2,$$ by \eqref{1.7} and
the Sobolev embedding $L^r\hookrightarrow\dH^{-3\alpha(r)}$.
\section{Proof of Proposition \ref{prop1.3}}\label{sec5}

The purpose of this section is to present the proof of Proposition
\ref{prop1.3}. Indeed it follows from   Proposition \ref{prop1.2}
that: for any $t\in[0,T]$,
\begin{align}\label{5.1}
\begin{split}
\cE(T)\cdot\Bigl(\int_0^t\|\nabla\p3{\V_+}& \|_{\htr}^2+\|\nabla\p3{\V_-}\|_{\htr}^2 dt'\Bigr)^{\frac r2}\\
& \leqslant\cE(T)\cdot\bigl(\|\Omega_0\|_{L^r}^r+\|j_0\|_{L^r}^r\bigr)+\Rmnum{3}_1(t)+\Rmnum{3}_2(t),
\end{split}
\end{align}
where
\begin{align*}
& \Rmnum{3}_1(t)\eqdefa\cE(T)\cdot\Bigl(\int_0^t(\|\u3\|_{\dH^{\frac 12+\frac 2p}}+\|\b3\|_{\dH^{\frac 12+\frac 2p}})\bigl(\|\odr\|_{L^2}^{2(2\alpha(r)+\frac 1p)}\|\nabla\odr\|_{L^2}^{2(1-\frac 1p)}\\
& \qquad\qquad\ +\|\omdr\|_{L^2}^{2(2\alpha(r)+\frac 1p)}\|\nabla\omdr\|_{L^2}^{2(1-\frac 1p)}\bigr)dt'\Bigr)^{\frac r2},\\
& \Rmnum{3}_2(t)\eqdefa\cE(T)\cdot\Bigl(\int_0^t(\|\u3\|_{\dH^{\frac 12+\frac 2p}}^2+\|\b3\|_{\dH^{\frac 12+\frac 2p}}^2)\bigl(\|\odr\|_{L^2}^2+\|\omdr\|_{L^2}^2\bigr)^{2(\alpha(r)+\frac 1p)}\\
&
\qquad\qquad\ \times\bigl(\|\nabla\odr\|_{L^2}^2+\|\nabla\omdr\|_{L^2}^2\bigr)^{1-\frac
2p}dt'\Bigr)^{\frac r2}.
\end{align*}
We emphasize that the constants in $\cE(t)$ may change from line to
line.

Applying H\"{o}lder's inequality gives \beno
\begin{split}
 |\Rmnum{3}_1(t)&|
\leqslant\cE(T)\Bigl(\int_0^t\|\nabla\odr\|_{L^2}^2dt'\Bigr)^{\frac
r2\cdot(1-\frac 1p)}\\
&\qquad\times\Bigl(\int_0^t(\|\u3\|^p_{\dH^{\frac 12+\frac
2p}}+\|\b3\|^p_{\dH^{\frac 12
+\frac 2p}})\|\odr\|_{L^2}^{2(1+2p\alpha(r))}dt'\Bigr)^{\frac r2\cdot\frac 1p}\\
&+\cE(T)\Bigl(\int_0^t\|\nabla\omdr\|_{L^2}^2dt'\Bigr)^{\frac r2\cdot(1-\frac 1p)}\\
&\qquad\times\Bigl(\int_0^t(\|\u3\|^p_{\dH^{\frac 12+\frac
2p}}+\|\b3\|^p_{\dH^{\frac 12+\frac 2p}})
\|\omdr\|_{L^2}^{2(1+2p\alpha(r))}dt'\Bigr)^{\frac r2\cdot\frac 1p},
\end{split} \eeno
then  Young's inequality yields \beq\label{5.2}\begin{split}
|&\Rmnum{3}_1(t)|
\leqslant\frac{r-1}{3r^2}\int_0^t\|\nabla\odr\|_{L^2}^2+\|\nabla\omdr\|_{L^2}^2dt'\\
&+\cE(T)\Bigl(\int_0^t(\|\u3\|^p_{\dH^{\frac 12+\frac
2p}}+\|\b3\|^p_{\dH^{\frac 12+\frac
2p}})\bigl(\|\odr\|_{L^2}+\|\omdr\|_{L^2}\bigr)^{2(1+2p\alpha(r))}dt'\Bigr)^{\frac{1}{1+2p\alpha(r)}}.
\end{split}
\eeq
 Similarly, we have
\begin{align}\label{5.3}
\begin{split}
|&\Rmnum{3}_2(t)|
\leqslant \frac{r-1}{3r^2}\int_0^t\|\nabla\odr\|_{L^2}^2+\|\nabla\omdr\|_{L^2}^2dt'\\
&+\cE(T)\Bigl(\int_0^t\bigl(\|\u3\|^p_{\dH^{\frac 12+\frac 2p}}
+\|\b3\|^p_{\dH^{\frac 12+\frac
2p}}\bigr)\bigl(\|\odr\|_{L^2}+\|\omdr\|_{L^2}\bigr)^{2(1+p\alpha(r))}dt'\Bigr)^{\frac{1}{1+p\alpha(r)}}.
\end{split}
\end{align}
For the last term, we get, by applying H\"{o}lder's inequality, that
\begin{align*}
 \Bigl(\int_0^t&\bigl(\|\u3\|^p_{\dH^{\frac 12+\frac 2p}}+\|\b3\|^p_{\dH^{\frac 12+\frac 2p}}\bigr)\bigl(\|\odr\|_{L^2}+\|\omdr\|_{L^2}\bigr)^{2(1+p\alpha(r))}dt'\Bigr)^{\frac{1}{1+p\alpha(r)}}\\
\leqslant& \Bigl(\int_0^t\bigl(\|\u3\|^p_{\dH^{\frac 12+\frac 2p}}+\|\b3\|^p_{\dH^{\frac 12+\frac 2p}}\bigr)dt'\Bigr)^{\frac{p\alpha(r)}{1+2p\alpha(r)}\cdot\frac{1}{1+p\alpha(r)}}\\
& \times\Bigl(\int_0^t\bigl(\|\u3\|^p_{\dH^{\frac 12+\frac
2p}}+\|\b3\|^p_{\dH^{\frac 12+\frac
2p}}\bigr)\bigl(\|\odr\|_{L^2}+\|\omdr\|_{L^2}\bigr)^{2(1+2p\alpha(r))}dt'\Bigr)^{\frac{1}{1+2p\alpha(r)}},
\end{align*}
and the definition of $\cE(T)$ implies
$$\cE(T)\cdot\Bigl(\int_0^t\bigl(\|\u3\|^p_{\dH^{\frac 12+\frac 2p}}+\|\b3\|^p_{\dH^{\frac 12+\frac 2p}}\bigr)dt'\Bigr)^{\frac{p\alpha(r)}{1+2p\alpha(r)}\cdot\frac{1}{1+p\alpha(r)}}\leqslant\cE(T).$$
Thus we deduce from \eqref{5.3} that
\begin{align}\label{5.4}
\begin{split}
|&\Rmnum{3}_2(t)|\leqslant \frac{r-1}{3r^2}\int_0^t\|\nabla\odr\|_{L^2}^2+\|\nabla\omdr\|_{L^2}^2dt'\\
& +\cE(T)\Bigl(\int_0^t\bigl(\|\u3\|^p_{\dH^{\frac 12+\frac
2p}}+\|\b3\|^p_{\dH^{\frac 12+\frac
2p}}\bigr)\bigl(\|\odr\|_{L^2}+\|\omdr\|_{L^2}\bigr)^{2(1+2p\alpha(r))}dt'\Bigr)^{\frac{1}{1+2p\alpha(r)}}.
\end{split}
\end{align}
Inserting \eqref{5.2} and \eqref{5.4} into \eqref{5.1} gives, for
any $t\in[0,T]$, \beno
\begin{split}
\cE(T)\Bigl(&\int_0^t\|\nabla\p3{\V_+}\|_{\htr}^2+\|\nabla\p3{\V_-}\|_{\htr}^2
dt'\Bigr)^{\frac r2}
\leqslant \frac{2(r-1)}{3r^2}\int_0^t\|\nabla\odr\|_{L^2}^2+\|\nabla\omdr\|_{L^2}^2dt'\\
& +\cE(T)\bigl(\|\Omega_0\|_{L^r}^r+\|j_0\|_{L^r}^r\bigr)+\cE(T)\Bigl(\int_0^t\bigl(\|\u3\|^p_{\dH^{\frac 12+\frac
2p}}+\|\b3\|^p_{\dH^{\frac 12+\frac
2p}}\bigr)\\
&\qquad\qquad\qquad\qquad\qquad\qquad\times\bigl(\|\odr\|_{L^2}+\|\omdr\|_{L^2}\bigr)^{2(1+2p\alpha(r))}dt'\Bigr)^{\frac{1}{1+2p\alpha(r)}}.
\end{split}
\eeno Then inserting the above inequality into the right hand side
of \eqref{1.8} gives
\begin{align*}
\frac 1r\bigl(&\|\odr(t)\|_{L^2}^2+\|\omdr(t)\|_{L^2}^2\bigr)+\frac{r-1}{r^2}\int_0^t\bigl(\|\nabla\odr\|_{L^2}^2+\|\nabla \omdr\|_{L^2}^2\bigr)dt'\\
\lesssim& \bigl(\frac 2r+1\bigr)\cE(T)\bigl(\|\Omega_0\|_{L^r}^r+\|j_0\|_{L^r}^r\bigr)\\
&+\cE(T)\Bigl(\int_0^t\bigl(\|\u3\|^p_{\dH^{\frac 12+\frac
2p}}+\|\b3\|^p_{\dH^{\frac 12+\frac
2p}}\bigr)\bigl(\|\odr\|_{L^2}+\|\omdr\|_{L^2}\bigr)^{2(1+2p\alpha(r))}dt'\Bigr)^{\frac{1}{1+2p\alpha(r)}}.
\end{align*}
Taking the power $1+2p\alpha(r)$ of this inequality and using the
elementary inequality
$$(a+b)^{\sigma}\thicksim a^{\sigma}+b^{\sigma},$$
for any positive index $\sigma$ and $a,\ b>0$, then we obtain for any $t\in[0,T]$,
\begin{align*}
 \|&\odr(t)\|_{L^2}^{2(1+2p\alpha(r))}+\|\omdr (t)\|_{L^2}^{2(1+2p\alpha(r))}\\
&\qquad+\Bigl(\int_0^t(\|\nabla\odr\|_{L^2}^2+\|\nabla \omdr\|_{L^2}^2)dt'\Bigr)^{1+2p\alpha(r)}\\
&\lesssim  \cE(T)\bigl(\|\Omega_0\|_{L^r}^{r(1+2p\alpha(r))}+\|j_0\|_{L^r}^{r(1+2p\alpha(r))}\bigr)\\
&\quad +\cE(T)\Bigl(\int_0^t\bigl(\|\u3\|^p_{\dH^{\frac 12+\frac
2p}}+\|\b3\|^p_{\dH^{\frac 12+\frac
2p}}\bigr)\bigl(\|\odr\|_{L^2}^{2(1+2p\alpha(r))}+\|\omdr\|_{L^2}^{2(1+2p\alpha(r))}\bigr)dt'\Bigr).
\end{align*}
Then Gronwall's inequality leads to \eqref{1.12}, which completes
the proof of the first part of Proposition \ref{prop1.3}.

Finally it follows from  Proposition \ref{prop1.2}, H\"{o}lder's
inequality and \eqref{1.12} that
\begin{align*}
& \bigl(\|\p3 {\V_+}(t)\|_{\htr}^2+\|\p3 {\V_-}(t)\|_{\htr}^2\bigr)+\int_0^t\bigl(\|\nabla\p3{\V_+}\|_{\htr}^2+\|\nabla \p3{\V_-}\|_{\htr}^2\bigr)dt'\\
\lesssim&\cE(t)\Bigl(\|\Omega_0\|_{L^r}^2+\|j_0\|_{L^r}^2\\
&+\bigl(\|\u3\|_{L_t^p(\dH^{\frac 12+\frac 2p})}+\|\b3\|_{L_t^p(\dH^{\frac 12+\frac 2p})}\bigr)\bigl(\|\odr\|_{L_t^{\infty}(L^2)}^{2(2\alpha(r)+\frac 1p)}+\|\omdr\|_{L_t^{\infty}(L^2)}^{2(2\alpha(r)+\frac 1p)}\bigr)\\
& \times\bigl(\|\nabla\odr\|_{L_t^2(L^2)}^{2(1-\frac 1p)}+\|\nabla\omdr\|_{L_t^2(L^2)}^{2(1-\frac 1p)}\bigr)+\bigl(\|\u3\|^2_{L_t^p(\dH^{\frac 12+\frac 2p})}+\|\b3\|^2_{L_t^p(\dH^{\frac 12+\frac 2p})}\bigr)\\
& \times\bigl(\|\odr\|_{L_t^{\infty}(L^2)}^{4(\alpha(r)+\frac
1p)}+\|\omdr\|_{L_t^{\infty}(L^2)}^{4(\alpha(r)+\frac 1p)}\bigr)
\bigl(\|\nabla\odr\|_{L_t^2(L^2)}^{2(1-\frac 2p)}+\|\nabla\omdr\|_{L_t^2(L^2)}^{2(1-\frac 2p)}\bigr)\Bigr)\\
\lesssim& \exp\bigl(C\cE(t)\bigr)\bigl(\|\Omega_0\|_{L^r}^2+\|j_0\|_{L^r}^2\bigr).
\end{align*}
This completes the proof of Proposition \ref{prop1.3}.

\section{Conclusion of the proof of Theorem \ref{thm1.3}}\label{sec6}
By Proposition \ref{prop1.3}, if we assume
\begin{equation}\label{6.1}
\int_0^{T^*}\|u^3\|^p_{\dH^{\frac 12+\frac 2p}}+\|b\|^p_{\dH^{\frac
12+\frac 2p}}dt'<\infty,
\end{equation}
we know that all the quantities in \eqref{6.2} are finite.
 We want to prove that all the above quantities prevent
the solution from blowing up. In order to do so, let us recall the
following theorem of anisotropic condition for blow up, which is a
generalization of Theorem 2.1 of \cite{CZ14} for the classical
Navier-Stokes system:

\begin{thm}\label{thm6.1}
(Proposition 4.1 of \cite{Yamazaki}) Let $u,b\in C([0,T^*[;\dH^{\frac
12}(\R^3))\bigcap L^2([0,T^*[;\dH^{\frac 32}(\R^3))$ solve the MHD
system \eqref{MHD}. If $T^*<\infty$, then for any
$p_{k,l}\in]1,\infty[,\ k,l\in\{1,2,3\}$, one has
\begin{equation}\label{6.3}
\sum_{k,\ell=1}^3 \int_0^{T^*}\bigl(\|\pa_\ell
u^k(t')\|_{\cB_{p_{k,l}}}^{p_{k,l}}+\|\pa_\ell
b^k(t')\|_{\cB_{p_{k,l}}}^{p_{k,l}}\bigr)dt'=\infty,
\end{equation}
where $\cB_p\eqdefa\dB_{\infty,\infty}^{-2+\frac 2p}$.
\end{thm}

Now let us present the proof of Theorem \ref{thm1.3}. Firstly, for
any $p\in]4,\infty[$,
\begin{align*}
\max_{1\leqslant
l\leqslant3}(\|\pa_\ell\u3\|_{\cB_p}+\|\pa_\ell\b3\|_{\cB_p})& \lesssim
\sup_{j\in\mathbb{Z}}2^{j(\frac 12+\frac 2p)}(\|\dj\u3\|_{L^2}+\|\dj\b3\|_{L^2})\\
& \lesssim\|\u3\|_{\dH^{\frac 12+\frac 2p}}+\|\b3\|_{\dH^{\frac
12+\frac 2p}}
\end{align*}
by Bernstein's inequality, which implies
\begin{equation}\label{6.4}
\max_{1\leqslant
l\leqslant3}\int_0^{T^*}\|\pa_\ell\u3\|_{\cB_p}^p+\|\pa_\ell\b3\|_{\cB_p}^p
dt' \lesssim\int_0^{T^*}\|\u3\|_{\dH^{\frac 12+\frac
2p}}^p+\|\b3\|_{\dH^{\frac 12+\frac 2p}}^p dt'\lesssim 1.
\end{equation}
Next, using Bernstein's inequality and the continuity of Riesz
transform in $L^p,\,\forall p\in]4,\infty[$, we have
\begin{align}\label{6.5}
\begin{split}
\int_0^{T^*}\|\nh\uh_{\div}\|_{\cB_p}^p+\|\nh\bh_{\div}\|_{\cB_p}^p
dt'
=& \int_0^{T^*}\|\nh\nh\Lh\p3\u3\|_{\cB_p}^p+\|\nh\nh\Lh\p3\b3\|_{\cB_p}^p dt'\\
\lesssim& \int_0^{T^*}\|\ph^2\Lh\u3\|_{\dH^{\frac 12+\frac 2p}}^p+\|\ph^2\Lh\b3\|_{\dH^{\frac 12+\frac 2p}}^p dt'\\
\thicksim & \int_0^{T^*}\|\u3\|_{\dH^{\frac 12+\frac
2p}}^p+\|\b3\|_{\dH^{\frac 12+\frac 2p}}^p dt' \lesssim 1.
\end{split}
\end{align}
The other components of the matrix $\nabla u$ and $\nabla b$ can be
estimated with norms which are not of scaling zero, namely norms
related to $\omega$ and d which have the scaling of $L^r$ norm as
showm in \eqref{6.2}. To proceed further, we first get for any
function a
\begin{align*}
\|\ddj a\|_{L^\infty}& \lesssim\sum_{k\leqslant j+1,l\leqslant j+1}2^k 2^{\frac l2}\|\dot{\Delta}_k^\h\dot{\Delta}_\ell^\v a\|_{L^2}\\
& \lesssim\|a\|_{\dH^{1-3\alpha(r)+\theta,-\theta}}\sum_{k\leqslant j+1,l\leqslant j+1}2^{k(3\alpha(r)-\theta)} 2^{l(\frac 12+\theta)}\\
& \lesssim2^{j(\frac
12+3\alpha(r))}\|a\|_{\dH^{1-3\alpha(r)+\theta,-\theta}},
\end{align*}
because $-(\frac 12+3\alpha(r))=-2+\frac{3}{r'}$, this leads to
\begin{equation}\label{6.6}
\|a\|_{\cB_q(r)}\lesssim\|a\|_{\dH^{1-3\alpha(r)+\theta,-\theta}},
\end{equation}
where $q(r)\eqdefa\frac{2r'}{3}$. As $r\in]\frac 32,2[$, q(r) is in
$]\frac 43,2[$. Applying mean inequality and triangle
inequality for the Besov norm, then \eqref{6.6}, H\"{o}lder's
inequality and \eqref{1.13}, we deduce that
\begin{align}\label{6.7}
\begin{split}
\int_0^{T^*}&\|\p3\uh_{\div}\|_{\cB_{q(r)}}^{q(r)}+\|\p3\bh_{\div}\|_{\cB_{q(r)}}^{q(r)}
dt'\\
\lesssim& \int_0^{T^*}\|\p3(\uh_{\div}+\bh_{\div})\|_{\cB_{q(r)}}^{q(r)}+\|\p3(\uh_{\div}-\bh_{\div})\|_{\cB_{q(r)}}^{q(r)} dt'\\
=& \int_0^{T^*}\|\nh\Lh\p3^2{\V_+}\|_{\cB_{q(r)}}^{q(r)}+\|\nh\Lh\p3^2{\V_-}\|_{\cB_{q(r)}}^{q(r)} dt'\\
\lesssim& \int_0^{T^*}\|\nh\Lh\p3^2{\V_+}\|_{\dH^{1-3\alpha(r)+\theta,-\theta}}^{q(r)}+\|\nh\Lh\p3^2{\V_-}\|_{\dH^{1-3\alpha(r)+\theta,-\theta}}^{q(r)} dt'\\
\lesssim&
(T^*)^{1-\frac{q(r)}{2}}\Bigl(\int_0^{T^*}\|\p3^2{\V_+}\|_{\htr}^2+\|\p3^2{\V_-}\|_{\htr}^2
dt'\Bigr)^{\frac{q(r)}{2}} \lesssim 1.
\end{split}
\end{align}
For the rest terms, firstly we get
\begin{align}\label{6.8}
\begin{split}
\|\nh\uh_{\curl}\|_{\cB_{q(r)}}+\|\nh\bh_{\curl}\|_{\cB_{q(r)}}
\lesssim& \|\nh(\uh_{\curl}+\bh_{\curl})\|_{\cB_{q(r)}}+\|\nh(\uh_{\curl}-\bh_{\curl})\|_{\cB_{q(r)}}\\
=& \|\nh\nh^{\bot}\Lh{\Gamma_+}\|_{\cB_{q(r)}}+\|\nh\nh^{\bot}\Lh{\Gamma_-}\|_{\cB_{q(r)}}\\
\lesssim& \|\ph^2\Lh{\Gamma_+}\|_{\dH^{1-3\alpha(r)}}+\|\ph^2\Lh{\Gamma_-}\|_{\dH^{1-3\alpha(r)}}\\
\lesssim& \|\nabla{\Gamma_+}\|_{\dH^{-3\alpha(r)}}+\|\nabla{\Gamma_-}\|_{\dH^{-3\alpha(r)}}\\
\lesssim& \|\nabla{\Gamma_+}\|_{L^r}+\|\nabla{\Gamma_-}\|_{L^r},
\end{split}
\end{align}
by \eqref{6.6}, continuity of Riesz transform in
$L^p,\ \forall p\in]1,\infty[$ and the Sobolev embedding
$L^r\hookrightarrow\dH^{-3\alpha(r)}$. Next, we use anisotropic
Bony's decomposition and Bernstein's inequality to get
\begin{align}\label{6.9}
\begin{split}
\|\ddj \p3(\uh_{\curl}\pm\bh_{\curl})\|_{L^\infty}& \lesssim\sum_{k\leqslant j+1,l\leqslant j+1}\|\ddj\dot{\Delta}_k^\h\dot{\Delta}_\ell^\v \p3\nh^{\bot}\Lh(\omega\pm d)\|_{L^{\infty}}\\
& \lesssim\sum_{k\leqslant j+1,l\leqslant j+1} 2^{k(\frac 2r-1)} 2^{\frac lr}\|\p3(\omega\pm d)\|_{L^r}\\
& \lesssim 2^{j(\frac 3r-1)}\|\p3(\omega\pm d)\|_{L^r}.
\end{split}
\end{align}
Recall $q(r)=\frac{2r'}{3}$, we find that $-(\frac
3r-1)=-2+\frac{2}{q(r)}$. Thus \eqref{6.9} actually leads to
\begin{equation}\label{6.10}
\|\p3(\uh_{\curl}\pm\bh_{\curl})\|_{\cB_{q(r)}}\lesssim\|\p3(\omega\pm
d)\|_{L^r}.
\end{equation}
Combining \eqref{6.8} and \eqref{6.10} gives
\begin{align}\label{6.11}
\begin{split}
 \|\nabla\uh_{\curl}\|_{\cB_{q(r)}}+\|\nabla\bh_{\curl}\|_{\cB_{q(r)}}
\lesssim& \|\nabla{\Gamma_+}\|_{L^r}+\|\nabla {\Gamma_-}\|_{L^r}\\
\lesssim& \|\nabla\odr\|_{L^2}\|\odr\|_{L^2}^{\frac
2r-1}+\|\nabla\omdr\|_{L^2}\|\omdr\|_{L^2}^{\frac 2r-1},
\end{split}
\end{align}
where we used \eqref{2.3} in the last step. Then combine
\eqref{6.11} with \eqref{1.12}, we get
\begin{align}\label{6.12}
\begin{split}
 \int_0^{T^*}\bigl(&\|\nabla\uh_{\curl}(t')\|_{\cB_{q(r)}}^{q(r)}+\|\nabla\bh_{\curl}(t')\|_{\cB_{q(r)}}^{q(r)}\bigr)
 dt'\\
\lesssim& (T^*)^{(1-\frac{q(r)}{2})}\|\odr\|_{L^{\infty}([0,T^*[;L^2)}^{q(r)\cdot(\frac 2r-1)}\bigl(\int_0^{T^*}\|\nabla\odr\|_{L^2}^2 dt'\bigr)^{\frac{q(r)}{2}}\\
&
\quad+(T^*)^{(1-\frac{q(r)}{2})}\|\omdr\|_{L^{\infty}([0,T^*[;L^2)}^{q(r)\cdot(\frac
2r-1)}\bigl(\int_0^{T^*}\|\nabla\omdr\|_{L^2}^2
dt'\bigr)^{\frac{q(r)}{2}} \lesssim 1.
\end{split}
\end{align}
Together with inequalities \eqref{6.4}, \eqref{6.5}, \eqref{6.7},
\eqref{6.12} and Theorem \ref{thm6.1}, we conclude the proof of
Theorem \ref{thm1.3}.

\noindent {\bf Acknowledgments.} This work was done when I was
visiting Morningside Center of the Academy of Mathematics and
Systems Sciences, CAS. I appreciate  the financial support from MCM.
I also would like to thank Professor Ping Zhang for guidance  and
careful reading of the preliminary version of this paper, and
Professor Zhifei Zhang for some enlightening comments.

\end{document}